\documentclass[lecture]{pcms-l}

\makeatletter

\newif\iffirstlecture\firstlecturefalse

\newcommand{\lectureseries}{\firstlecturetrue
              \secdef\@lectureseries\@slectureseries} 

\newcommand{\@lectureseries}[2][default]{\chapter*{#2}%
              \gdef\thelectureseries{#1}} 

\newcommand{\@slectureseries}[1]{\chapter*{#1}}

\renewcommand{\auth}{\secdef\@auth\@sauth}

\newcommand{\@auth}[2][default]{\vspace{-1pc}{\raggedleft
        \Large\bf\noindent
        #2\endgraf
        \vspace*{2pc}
        }
        \def\@author{#1}%
}

\newcommand{\@sauth}[1]{\vspace{-1pc}{\raggedleft
        \Large\bf\noindent
        #1\endgraf
        \vspace*{2pc}
        }
        \def\@author{#1}%
}

\def\lecture#1{\global\Lecturetrue\global\Monographfalse
\iffirstlecture\else\chapter*{}\fi\firstlecturefalse
  \global\let\sectionmark\@gobble 
  \addtocounter{lecture}1\relax
  \refstepcounter{chapter}%
\gdef\thelecturename{#1\unskip}
  {\Large\bfseries
   \raggedleft
   \@xp\uppercase\@xp{\thelecturelabel} {\LARGE\thelecturenum}\\
   \vspace*{3pt}%
   \thelecturename
   \endgraf}%
  \let\@secnumber=\thelecturenum
  \@xp\lecturemark\@xp{\thelecturename}%
  \addcontentsline{toc}{chapter}{%
    \thelecturelabel\ \thelecturenum.\ \thelecturename}%
  \vspace{34\p@}\noindent}
  
\def\lecture{\global\Lecturetrue\global\Monographfalse
\iffirstlecture\else\chapter*{}\fi%
  \global\let\sectionmark\@gobble 
\secdef\@lecture\@slecture}

\def\@lecture[#1]#2{%
  \addtocounter{lecture}1\relax
  \refstepcounter{chapter}%
\gdef\thelecturename{#1\unskip}\firstlecturefalse
  {\Large\bfseries
   \raggedleft
   \@xp\uppercase\@xp{\thelecturelabel} {\LARGE\thelecturenum}\\
   \vspace*{3pt}%
    #2\unskip
   \endgraf}%
  \let\@secnumber=\thelecturenum
  \@xp\lecturemark\@xp{\thelecturename}%
  \addcontentsline{toc}{chapter}{%
    \thelecturelabel\ \thelecturenum.\ #2}%
  \vspace{34\p@}\noindent}
  
\def\slecturerunhead#1#2#3{%
    \let\@tempa\chaptername
    \uppercasenonmath{\@tempa}%
    \def\@tempb{#3\unskip}%
    \uppercasenonmath{\@tempb}%
    {\normalfont\@tempb}
    }
\def\slecturemark{
    \@secmark\markright\slecturerunhead\chaptername}%

\def\@slecture#1{%
\iffirstlecture
\gdef\thelecturename{#1\unskip}\firstlecturefalse
  {\Large\bfseries
\noindent\thelecturename
   \endgraf}%
  \let\@secnumber=\thelecturenum
  \@xp\slecturemark\@xp{\thelecturename}%
  \addcontentsline{toc}{chapter}{%
    \thelecturename}%
 \vspace{-6\p@}\noindent
\else
\gdef\thelecturename{#1\unskip}\firstlecturefalse
  {\Large\bfseries
   \raggedleft
   \@xp\uppercase\@xp{\thelecturename}
   \endgraf}%
  \let\@secnumber=\thelecturenum
  \@xp\slecturemark\@xp{\thelecturename}%
  \addcontentsline{toc}{chapter}{%
    \thelecturename}%
  \vspace{34\p@}\noindent
\fi}

\ifLecture
  \def\chapterrunhead#1#2#3{%
    \let\@tempa\@author
    \uppercasenonmath{\@tempa}%
    \uppercasenonmath{\thelectureseries}%
    \textmd{\@tempa, \thelectureseries}
    }
  \def\lecturerunhead#1#2#3{%
    \let\@tempa\chaptername
    \uppercasenonmath{\@tempa}%
    \def\@tempb{#3\unskip}%
    \uppercasenonmath{\@tempb}%
    \textmd{\@tempa\ #2. \@tempb}
    }
\else
  \let\chapterrunhead\partrunhead
\fi

\newif\ifBibliographyIsASection\BibliographyIsASectionfalse

  \def\bibliomark{
    \@secmark\markright\bibliorunhead\chaptername}%

  \def\bibliorunhead#1#2#3{%
    \let\@tempa\chaptername
    \uppercasenonmath{\@tempa}%
    \def\@tempb{#3\unskip}%
    \uppercasenonmath{\@tempb}%
    \textmd{\@tempb}
    }

\def\thebibliography#1{%
  \ifBibliographyIsASection
    \section*\refname
    \if@backmatter
      \markboth{\refname}{\refname}%
    \fi
  \else
\chapter*{}
  {\Large\bfseries
   \raggedleft
   \@xp\uppercase\@xp{\bibname} \\
   \endgraf}%
  \let\@secnumber=\thelecturenum
  \@xp\bibliomark\@xp{\bibname}%
  \addcontentsline{toc}{chapter}{%
    \bibname}%
  \vspace{34\p@}\noindent
  \fi
  \normalsize\labelsep .5em\relax
  \list{\@arabic\c@enumi.}{\settowidth\labelwidth{\@biblabel{#1}}%
  \leftmargin\labelwidth
  \advance\leftmargin\labelsep
	\usecounter{enumi}}\sloppy
  \clubpenalty9999 \widowpenalty\clubpenalty  \sfcode`\.\@m}

  \def\indexmark{
    \@secmark\markright\indexrunhead\chaptername}%

  \def\indexrunhead#1#2#3{%
    \let\@tempa\chaptername
    \uppercasenonmath{\@tempa}%
    \def\@tempb{#3\unskip}%
    \uppercasenonmath{\@tempb}%
    \textmd{\@tempb}
    }

\ifLecture
\def\theindex{\cleardoublepage
\@restonecoltrue\if@twocolumn\@restonecolfalse\fi
\columnseprule \z@ \columnsep 35pt
\def\indexchap{\@startsection
		{chapter}{1}{\z@}{8pc}{34pt}%
		{\raggedleft
		\Large\bfseries}}%
 \twocolumn[\indexchap[{\indexname}]{\@xp\uppercase\@xp{\indexname}}]
  \@xp\indexmark\@xp{\indexname}%
	\thispagestyle{plain}\let\item\@idxitem\parindent\z@
	 \footnotesize\parskip\z@ plus .3pt\relax\let\item\@idxitem}
\fi

\def\@makefntext{\noindent\@makefnmark}

\def\setaddress{%
  {\let\@makefnmark\relax  \let\@thefnmark\relax
        \nobreak
        \addressnum@=\z@
        \loop\ifnum\addressnum@<\addresscount@\advance\addressnum@\@ne
           \footnote{$^{\hbox{\tiny\number\addressnum@}}$%
           \csname @address\number\addressnum@\endcsname
           \csname @curraddr\number\addressnum@\endcsname
           \csname @email\number\addressnum@\endcsname}\repeat
  \ifx\@empty\@date\else \@footnotetext{\@setdate}\fi
  \ifx\@empty\@subjclass\else \@footnotetext{\@setsubjclass}\fi
  \ifx\@empty\@keywords\else \@footnotetext{\@setkeywords}\fi
  \ifx\@empty\thankses\else \@footnotetext{%
    \def\par{\let\par\@par}\@setthanks}\fi
    }%
  \@setcopyright
}

\def\@tmpevenhead{\relax}

\def\cleardoublepage{\clearpage\if@twoside \ifodd\c@page\else
    \let\@tmpevenhead\@evenhead \let\@evenhead\relax\hbox{}\eject 
    \let\@evenhead\@tmpevenhead\if@twocolumn\hbox{}\newpage\fi\fi\fi}

\def\@setcopyright{%
  \let\copyrightyear\currentyear             
  \insert\copyins{\hsize\textwidth
    \parfillskip\z@ \leftskip\z@\@plus.9\textwidth
    \fontsize{6}{7\p@}\normalfont\upshape
    \everypar{}%
    \vskip-\skip\copyins \nointerlineskip
    \noindent\vrule\@width\z@\@height\skip\copyins
    \copyright\copyrightyear\ American Mathematical Society\par
    \kern\z@}%
}

\renewcommand{\@auth}[2][default]{{\raggedleft
        \begingroup
  \fontsize{\@xivpt}{18}\bfseries
  #2\par \endgroup
        \vspace*{2pc}
        }
        \def\@author{#1}%
}

\renewcommand{\@sauth}[1]{{\raggedleft
        \begingroup
  \fontsize{\@xivpt}{18}\bfseries
  #1\par \endgroup
        \vspace*{2pc}
        }
        \def\@author{#1}%
}

\def\@lecture[#1]#2{%
  \addtocounter{lecture}1\relax
  \refstepcounter{chapter}%
\gdef\thelecturename{#1\unskip}\firstlecturefalse
  {\Large\bfseries
   \raggedleft
   \@xp\uppercase\@xp{\thelecturelabel} {\LARGE\thelecturenum}\\
   \vspace*{3pt}%
    #2\unskip
   \endgraf}%
  \let\@secnumber=\thelecturenum
  \@xp\lecturemark\@xp{\thelecturename}%
  \addcontentsline{toc}{chapter}{%
    \thelecturelabel\ \thelecturenum.\ #2}%
  \vspace{10\p@}\noindent}
  
\def\@slecture#1{%
\iffirstlecture
\gdef\thelecturename{#1\unskip}\firstlecturefalse
  {\Large\bfseries
\noindent\thelecturename
   \endgraf}%
  \let\@secnumber=\thelecturenum
  \@xp\slecturemark\@xp{\thelecturename}%
  \addcontentsline{toc}{chapter}{%
    \thelecturename}%
 \vspace{-6\p@}\noindent
\else
\gdef\thelecturename{#1\unskip}\firstlecturefalse
  {\Large\bfseries
   \raggedleft
   \@xp\uppercase\@xp{\thelecturename}
   \endgraf}%
  \let\@secnumber=\thelecturenum
  \@xp\slecturemark\@xp{\thelecturename}%
  \addcontentsline{toc}{chapter}{%
    \thelecturename}%
  \vspace{10\p@}\noindent
\fi}

\makeatother

\usepackage{url,graphicx}
\usepackage[dvips]{epsfig}
\usepackage{latexsym,color}
\usepackage{amscd,amssymb,verbatim}

\newcommand{\br}{{\mathbb R}}
\newcommand{\rn}{{\mathbb R}^n}

\newcommand{\dz}{{\mathbb Z}}
\newcommand{\zz}{{{\mathbb Z}_2}}

\newcommand{\cc}{{\mathcal C}}
\newcommand{\cd}{{\mathcal D}}
\newcommand{\cf}{{\mathcal F}}

\newcommand{\ch}{{\mathcal H}}
\newcommand{\cm}{{\mathcal M}}
\newcommand{\cn}{{\mathcal N}}
\newcommand{\co}{{\mathcal O}}

\newcommand{\car}{{\mathcal R}}
\newcommand{\cs}{{\mathcal S}}
\newcommand{\ct}{{\mathcal T}}

\newcommand{\tth}{{\tt h}}

\newcommand{\aut}{\text{\rm Aut}\,}
\newcommand{\bd}{{\text{\rm Bd}\,}}
\newcommand{\bip}{{\text{\rm Bip}\,}}
\newcommand{\be}{\begin{enumerate}}
\newcommand{\bo}{\partial}
\newcommand{\bu}{\bullet}

\newcommand{\chom}{\text{\tt Hom}_{\,0}}
\newcommand{\coind}{\text{\rm Coind}\,}

\newcommand{\conn}{\mbox{\rm conn}}
\newcommand{\cmp}{\complement}
\newcommand{\da}{\Delta}
\newcommand{\ee}{\end{enumerate}}
\newcommand{\graphs}{{{\bf Graphs}}}
\newcommand{\hra}{\hookrightarrow}
\newcommand{\id}{\text{\rm id}}
\newcommand{\im}{\text{\rm im}\,}
\newcommand{\ind}{\text{\rm Ind}\,}
\newcommand{\join}[1]{\ast_{#1}}

\newcommand{\link}{\text{\rm lk}}
\newcommand{\lra}{\longrightarrow}
\newcommand{\mgraphs}{\text{\bf Graphs}}
\newcommand{\graphsb}{{\bf Graphs}$_{\text{\bf p}}$}

\newcommand{\nin}{\noindent}

\newcommand{\ovr}[1]{\overline{#1}}
\newcommand{\pr}{\noindent{\bf Proof. }}
\newcommand{\ra}{\rightarrow}
\newcommand{\rp}{{\mathbb R\mathbb P}}
\newcommand{\rk}{\text{\rm rk}\,}

\newcommand{\sm}{\setminus}
\newcommand{\stab}{\text{\rm stab}}
\newcommand{\supp}{\text{\rm supp}\,}

\newcommand{\sw}{\varpi_1}

\newcommand{\thom}{\text{\tt Hom}\,}
\newcommand{\thomp}{\text{\tt Hom}_+}
\newcommand{\ti}{\tilde}

\newcommand{\vt}{\vartheta}
\newcommand{\wind}{{\text{\rm wind}\,}}
\newcommand{\wti}{\widetilde}

\numberwithin{section}{lecture}

\newtheorem{thm}{Theorem}[section]

\newtheorem{df}[thm]{Definition}

\newtheorem{crl}[thm]{Corollary}
\newtheorem{prop}[thm]{Proposition}
\newtheorem{conj}[thm]{Conjecture}
\newtheorem{rem}[thm]{Remark}

\numberwithin{equation}{section}
\numberwithin{figure}{section}
\numberwithin{table}{section}

\begin{document}



\part*{Chromatic numbers, morphism complexes, and 
Stiefel-Whitney characteristic classes.}
\pauth{Dmitry N.\ Kozlov}
\tableofcontents

\mainmatter
\setcounter{page}{1}

\LogoOn

\lectureseries[Morphism complexes, and 
Stiefel-Whitney classes]
{Chromatic numbers, morphism complexes, and 
Stiefel-Whitney characteristic classes. }

\auth[D.\ N.\ Kozlov]{Dmitry N. Kozlov} 

\address{Institute of Theoretical Computer Science / Department of Mathematics, Eidgen\"ossische Technische Hochschule - Z\"urich, CH-8006 Z\"urich, Switzerland}
\email{dkozlov@inf.ethz.ch}
\thanks{The author would like to thank the Swiss National Science Foundation and Mathematical Science Research Institute, Berkeley for the 
generous support.}



\setaddress

\date\today

\nin {\sc Preamble.} Combinatorics, in particular graph theory, has a~rich history of being a~domain of successful applications of tools from other areas of mathematics, including topological methods. Here, we survey the study of the $\thom$-complexes, and the ways these can be used to obtain lower bounds for the chromatic numbers of graphs, presented in a~recent series of papers \cite{BK03a,BK03b,BK03c,CK1,CK2,K4,K5}. 

The structural theory is developed and put in the historical context,
culminating in the proof of the Lov\'asz Conjecture, which can be
stated as follows: 
\begin{quotation}{\it For a~graph $G$, such that the complex
$\thom(C_{2r+1},G)$ is $k$-connected for some $r,k\in\dz$, $r\geq 1$,
$k\geq -1$, we have $\chi(G)\geq k+4$.}\end{quotation}

Beyond the, more customary in this area, cohomology groups, the
algebro-topological concepts involved are spectral sequences and
Stiefel-Whitney characteristic classes.  Complete proofs are included
for all the new results appearing in this survey for the first time.


\newpage

\lecture{Introduction.}

\section{The chromatic number of a graph.}

\subsection{The definition and applications.}
\label{ss1.1} 


\nin Unless stated otherwise, all graphs are undirected, loops 
are allowed, whereas multiple edges are not. We shall occasionally
stress these conventions, to avoid the possibility of
misunderstanding.

For a graph $G$, $V(G)$ denotes the set of its vertices, and $E(G)$
denotes the set of its edges. If convenient, we think of $E(G)$ as
a~$\zz$-invariant subset of $V(G)\times V(G)$, where $\zz$ acts on
$V(G)\times V(G)$ by switching the coordinates: $(x,y)\mapsto(y,x)$.
Under this convention, a~looped vertex $x$ is encoded by the diagonal
element $(x,x)$, while the edge from $x$ to $y$ (for $x\neq y$) is
encoded by the pair $(x,y),(y,x)\in V(G)\times V(G)$. For example, the
edge set of the graph with 2 vertices connected by an edge, were the
first vertex is looped, and the second one is not, is encoded by the
set $\{(1,1),(1,2),(2,1)\}$.

\begin{df}\label{df:vrtcl}
Let $G$ be a graph. A {\bf vertex-coloring} of $G$ is 
a~set map $c:V(G)\rightarrow S$ such that $(x,y)\in E(G)$
implies $c(x)\neq c(y)$.
\end{df}

Clearly, a vertex coloring exists if and only if $G$ has no loops.

\begin{df}\label{df:chrnm}
The {\bf chromatic number} of~$G$, $\chi(G)$, is the minimal 
cardinality of a~finite set $S$, such that there exists a~vertex-coloring
$c:V(G)\rightarrow S$.
\end{df}

If no such finite set $S$ exists, for example, if $G$ has loops, we
use the convention $\chi(G)=\infty$.

\begin{figure}[hbt]
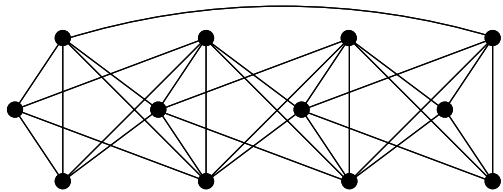

\begin{center}
  \begin{picture}(0,0)%
    \includegraphics{gr1.pstex}%
  \end{picture}%
  \input{gr1.pstex_t}%
 
\end{center}
\caption{A graph with chromatic number 4, which does not contain $K_4$ as an induced subgraph.}
\label{fig:chr1}
\end{figure}

The literature devoted to the applications of computing the chromatic number of a~graph is very extensive. Two of the basic applications are the {\it frequency assignment problem} and the {\it task scheduling problem}. 

The first one concerns a~collection of transmitters, with certain pairs of transmitters required to have different frequencies (e.g., because they are too close). Clearly, the minimal number of frequencies required for such an~assignment is precisely the chromatic number of the graph, whose vertices correspond to the transmitters, with two vertices connected by an~edge if and only if the corresponding transmitters are requested to have different frequencies. 

The second problem concerns a~collection of tasks which need to be performed. Each task has to be performed exactly once, and the tasks are to be performed in regularly allocated slots (e.g., hours). The only constraint is that certain tasks cannot be performed simultaneously. Again, the minimal number of slots required for the task scheduling is equal to the chromatic number of the graph, whose vertices correspond to tasks, with two vertices connected by an edge if and only if the corresponding tasks cannot be performed simultaneously.

\subsection{The Hadwiger Conjecture.}
\label{ss1.2} 


\nin The question of computing $\chi(G)$ has a long history. In 1852, F.\ Guthrie, \cite{Gut}, asked whether it is true that any planar map of connected countries can be colored with 4 colors, so that every pair of countries, which share a (non-point) boundary segment, receive different colors. 
The first time this question appeared in print was in a~paper by
Cayley, \cite{Cay}, and it became known as the Four-Color Problem, one
of the most famous questions in graph theory, as well as a~popular
brain-teaser. There is very extensive literature on the subject, see
e.g., \cite{Har,KS,MSTY,Ore1,Th1}.

The apparently first proof, offered by A.\ Kempe in 1880, \cite{Kem},
turned out to be false, as did many later ones. The flaw was noticed
in 1890 by P.\ Heawood, \cite{Hea}, who also proved the weaker
Five-Color Theorem. After several important contributions, most
notably by G.\ Birkhoff, and H.\ Heesch, \cite{Hee}, the latter
reduced the Four-Color Theorem to the analysis of the large, but
finite set of "unavoidable" configurations, the original conjecture
has been proved 1976 by Appel \& Haken, using computer computations,
see \cite{AH1} for the original announcement and \cite{AH} for last
reprint. A~new, shorter and more structural proof (though still
relying on computers) has been obtained in 1997, see \cite{RSST}. The
usual way to formulate this theorem is to dualize the map to obtain
a~planar graph, coloring vertices instead of the countries.

\begin{thm} \label{thm:4color} {\bf (The Four-Color Theorem).} {\rm (Appel \& Haken, \cite{AH}; revised proof by Robertson, Sanders, Seymour~\& Thomas, \cite{RSST}).}

\nin Every planar graph is four colorable.
\end{thm}
 
In 1943, Hadwiger, \cite{Had}, stated a conjecture closely related to the Four-Color Theorem. Recall that a~graph $H$ is called a~{\it~minor} of another graph $G$, if $H$ can be obtained from a~subgraph of $G$ by a~sequence of edge-contractions. Let $K_n$ denote an unlooped complete graph on $n$ vertices, that is, $V(K_n)=[n]$, $E(K_n)=\{(x,y)\,|\,x,y\in [n], x\neq y\}$.

 \begin{conj} \label{conj:hadw}
{\bf (Hadwiger Conjecture).} $\,$ 

\nin For every positive integer $t$, if a~graph has no $K_{t+1}$ minor, then it has a~$t$-coloring, in other words, every graph $G$ has $K_{\chi(G)}$ as its minor.
\end{conj}

The Hadwiger conjecture is proved for $\chi(G)\leq 5$. Indeed, it is trivial for $\chi(G)=1$, as $K_1$ is a~minor of any graph. For $\chi(G)=2$ it just says that $K_2$ is a~minor of an~arbitrary graph containing an~edge. If $\chi(G)=3$, then $G$ contains an odd cycle, in particular it has $K_3$ as a~minor. The case $\chi(G)=4$ is reasonably easy, and was shown by Hadwiger, \cite{Had}, and Dirac, \cite{Dir}. Finally, it was shown in 1937 by Wagner, \cite{Wag}, that the case $\chi(G)=5$ of the Hadwiger Conjecture is equivalent to the Four-Color Theorem.

\subsection{The complexity of computing the chromatic number.}
\label{ss1.3} 


\nin The problem of computing the chromatic number of a~graph is
NP-complete, implying that the worst-case performance of any algorithm
is, most likely, exponential in the number of vertices. Stronger,
already the, seemingly much more special, problem of deciding whether
a~given planar graph is 3-colorable is NP-complete, see e.g.,
\cite{GJ2}.  Recently, it has been shown that even coloring
a~3-colorable graph with 4 colors is NP-complete, see~\cite{KLS}.

We note that deciding whether $\chi(G)=2$ (i.e., whether the graph is
bipartite) is computationally much easier. The plain depth-first
search yields $O(|V(G)|+|E(G)|)$ performance time. Another good news
is that one can 4-color a~planar graph in polynomial time: quartic
time was obtained in \cite{AH}, and later improved to quadratic time,
see~\cite{RSST}.

The situation is not getting much better if we switch to considering
approximations. For example, it was shown by Garey \& Johnson,
\cite{GJ}, that if a~polynomial time approximate algorithm for graph
coloring exists (in the precise formulation, meaning that the output
of the algorithm does not differ by more than a~constant factor, which
is smaller than 2, from the actual value of the chromatic number),
then there exists a~polynomial time algorithm for graph coloring,
which, of course, is not very likely.

Much of the same can be said about computing lower bounds. Even the
most trivial lower bound for the chromatic number, given by the clique
number, is not good, since computing the clique number of a~graph is
also an~NP-complete problem. For fixed clique size, the lower bound
based on clique number is polynomially computable, but is not very
interesting.

The original Lov\'asz bound, by  virtue of being based on computing
the connectivity of a~simplicial complex, also has a~very high
computational complexity, since determining the triviality of the
homotopy groups is an~extremely hard problem, even in low dimensions.

It is then a positive and welcome surprise, that our bounds, based on
the Stiefel-Whitney characteristic classes are both nontrivial and
polynomially computable; here we fix the test graph and the tested
dimension and consider the computational complexity with respect to
the number of vertices of the graph which is being tested. The crucial
difference is that, as opposed to homotopy, the cohomology groups (and
the functorial invariants contained therein) may be computed by means
of simple linear algebra.

\section{The category of graphs.}

\subsection{Graph homomorphisms and the chromatic number.} 
\label{ss2.1} 


\nin The following notion is the gist in recasting the various coloring
questions in the functorial language.

\begin{df}
For two graphs $T$ and $G$, a {\bf graph homomorphism} from $T$ to $G$
is a~map $\varphi:V(T)\rightarrow V(G)$, such that if $x,y\in V(T)$
are connected by an edge, then $\varphi(x)$ and $\varphi(y)$ are also
connected by an edge.
\end{df}
In other words, the map of the vertex sets $\varphi:V(T)\rightarrow
V(G)$ induces the product map $\varphi\times\varphi:V(T)\times
V(T)\rightarrow V(G)\times V(G)$, and the condition for the set map
$\varphi$ to be a~graph homomorphism translates into
\[(\varphi\times\varphi)(E(T))\subseteq E(G).\] Expressed verbally:
{\it edges map to edges}.

The study of graph homomorphisms is a~classical and well-developed subject within combinatorics. The interested reader may want to consult the textbooks \cite{GoR} and \cite{HN04}.

Clearly, for any two positive integers $m$ and $n$, a~graph homomorphism $\varphi:K_m\ra K_n$ exists if and only if $m\leq n$. More generally, we can now restate Definition~\ref{df:chrnm} in the language of graph homomorphisms.

\begin{df}\label{df:chrnm2}
The {\bf chromatic number} of~$G$, $\chi(G)$, is the minimal positive
integer $n$, such that there exists a~graph homomorphism
$\varphi:G\rightarrow K_n$.
\end{df}

In this sense, the problem of vertex-colorings and computing chromatic
numbers corresponds to choosing a particular family of graphs, namely
unlooped complete graphs, fixing a~valuation on this family, here we
are mapping $K_n$ to $n$, and then searching for a graph homomorphism
from a given graph to the chosen family, which would minimize the
fixed valuation. Using the intuition from statistical mechanics we
call such a family of graphs {\it state graphs}.

A natural question arises: are there any other choices of families of
state graphs and valuations which correspond to other natural and
well-studied classes of graph problems. The answer is yes, and we
shall describe two examples in the following subsections.

\subsection{The fractional chromatic number.} \label{ss2.2} 


\nin First, we define an important family of graphs.

\begin{df}\label{df:kngr}
Let $n,k$ be positive integers, $n\geq 2k$. The {\bf Kneser graph}
$K_{n,k}$ is defined to be the graph whose set of vertices is the set
of all $k$-subsets of $[n]$, and the set of edges is the set of all
pairs of disjoint $k$-subsets.
\end{df}

\nin{\bf Examples:} 
\begin{itemize}
\item $K_{2k,k}$ is a matching on $\binom{2k}{k}$ vertices; 
\item $K_{n,1}$ is the unlooped complete graph $K_n$;
\item $K_{5,2}$ is the Petersen graph.
\end{itemize}

We can now define the fractional chromatic number by means of graph
homomorphisms.

\begin{df} \label{df:fracnm}
Let $G$ be a graph. The {\bf fractional chromatic number} of~$G$,
$\chi_f(G)$, is defined by 
\[\chi_f(G)=\inf_{(n,k)}\frac{n}{k}\,\,,
\]
where the infimum is taken over all pairs $(n,k)$ such that
there exists a graph homomorphism from $G$ to $K_{n,k}$.
\end{df}

Here, the state graphs are the Kneser graphs, $\{K_{n,k}\}_{n\geq
2k}$, and the chosen valuation on this family is $K_{n,k}\mapsto
n/k$.

\subsection{The circular chromatic number.} \label{ss2.3} 


\nin Again, we start by defining the appropriate family of graphs.

\begin{df}\label{df:rr}
Let $r$ be a real number, $r\geq 2$. $R_r$ is defined to be the graph
whose set of vertices is the set of unit vectors in the plane pointing
from the origin, and two vertices $x$ and $y$ are connected by an edge
if and only if $2\pi/r\leq\alpha$, where $\alpha$ is the sharper of
the two angles between $x$ and $y$ (or $\pi$ if these two angles are
equal).
\end{df}

Note that both the number of vertices and valencies of the vertices
(if $r>2$) are infinite.
 
\begin{df}\label{df:circnm}
Let $G$ be a graph. The {\bf circular chromatic number} of~$G$ is
$\chi_c(G)=\inf r$, where the infimum is taken over all positive reals
$r$, such that there exists a graph homomorphism from $G$ to $R_r$.
\end{df} 

In other words, the family of the state graphs is $\{R_r\}_{r\geq 2}$,
and the chosen valuation is $R_r\mapsto r$.

It is possible to define $\chi_c(G)$ by using only finite state graphs.

\begin{df}\label{df:rgr}
Let $n,k$ be positive integers, $n\geq 2k$. $R_{n,k}$ is defined to be
the graph whose set of vertices is $[n]$, and two vertices $x,y\in
[n]$ are connected by an edge if and only if 
\[k\leq|x-y|\leq n-k.\]
\end{df}

\nin{\bf Examples:} 
\begin{itemize}
\item $R_{2k,k}$ is a complete matching on $2k$ vertices; 
\item $R_{2k+1,k}$ is a cycle with $2k+1$ vertices;
\item $R_{n,1}=K_{n,1}=K_n$;
\item $R_{n,2}$ is the unlooped complement of a cycle with $n$ vertices. 
\end{itemize}

The equivalent definition of $\chi_c(G)$ in terms of finite state
graphs is also functorial.

\begin{prop}\label{df:circnm2}
Let $G$ be a graph. We have the equality
\[\chi_c(G)=\inf_{(n,k)}\frac{n}{k}\,\,,
\]
where the infimum is taken over all pairs $(n,k)$ such that there
exists a graph homomorphism from $G$ to $R_{n,k}$.
\end{prop}

Here, the state graphs are $\{R_{n,k}\}_{n\geq 2k}$, and the chosen
valuation on this family is again $R_{n,k}\mapsto n/k$. 

We remark, that for any graph $G$ we have
$$\chi(G)-1<\chi_c(G)\leq\chi(G).$$ For more information on the
circular chromatic number of a~graph, see~\cite{Vi,Zhu}.

\subsection{The category $\graphs$.} 
\label{ss2.4} 


\nin As we have seen so far, graph homomorphisms are invaluable in
formulating various coloring problems. The usual framework for
studying a~set of mathematical objects and the set of maps between
them is that of {\em a~category}. Before we get that, we need to check
a~few properties.

\nin (1) Let $T,G,H$ be three graphs, and let $\varphi:T\ra G$ and
$\psi:G\ra H$ be two graph homomorphisms. The {\em composition } of the set maps $\psi\circ\varphi$ is again a~graph homomorphism from $T$ to $H$, as $(x,y)\in E(T)$ implies $(\varphi(x),\varphi(y))\in E(G)$, which further implies $(\psi(\varphi(x)),\psi(\varphi(y)))\in E(H)$.

\nin (2) The composition of set maps (and hence of graph
homomorphisms) is {\em associative}.

\nin (3) For any graph $G$, the {\it identity} map $\id:V(G)\ra V(G)$
is a~graph homomorphism.

Now we are ready to put it all together into one structure.

\begin{df}\label{df:grcat}
{\graphs}  is the category defined as follows:
\begin{itemize}
\item the objects of $\mgraphs$ are all graphs;
\item the morphisms $\cm(G,H)$ for two objects $G,H\in\co(\mgraphs)$ 
are all graph homomorphisms from $G$ to $H$.
\end{itemize}
\end{df}

For a graph $G$, let $G^o$ be the {\it looping} of $G$, i.e.,
$V(G^o)=V(G)$, $E(G^o)=E(G)\cup\{(x,x)\,|\,x\in V(G)\}$.
Then, $K_1^o$ is a~graph consisting of one vertex and one loop, it is
the {\it terminal object} of {\bf~Graphs}. The empty graph is the
{\it~initial object} of {\bf~Graphs}.

As a useful variation we also consider the category \graphsb (where {\bf p} stand for "proper"), whose objects are all graphs, and whose morphisms are all proper graph homomorphisms. We call a~graph homomorphism $\varphi:T\ra G$ {\it proper} if $|\varphi^{-1}(g)|$ is finite for all $g\in V(G)$.

\begin{prop} \label{pr:grprcopr}
The direct product of graphs (see Definition~\ref{df:dirprod}) is a~categorical {\it~product}, while the disjoint union of graphs is a~categorical {\it~coproduct} in {\bf~Graphs}.
\end{prop}

Even more generally, we have the following property.
\begin{prop} \label{pr:}
{\bf Graphs} has all finite limits and colimits.
\end{prop}

A surprising result of Welzl, \cite{Welzl}, shows that this category is in a~certain sense dense.

\begin{thm} \label{thm:welzl} {\rm (Welzl Theorem).} $\,$

\nin Let $T$ and $G$ be two arbitrary finite graphs, such that $\chi(T)\geq 3$, and there exists a~graph homomorphism from $T$ to $G$, but there is no graph homomorphism from $G$ to~$T$. Then, there exists a~graph~$H$, such that there exist graph homomorphisms from $T$ to $H$, and from $H$ to $G$, but there are no graph homomorphisms from $H$ to $T$, or from $G$ to~$H$.
\end{thm}

In the setting of this category, we can think of the following generalization of various coloring problems: given a~category $\cc$, and a~subcategory $\wti\cc$, determine the set of morphisms from a~given object $A$ to the objects in $\wti\cc$. In other words, we need to study obstructions to the existence of morphisms between certain objects.

\subsection{Test objects.} 
\label{ss2.5} 


\nin Due to the lack of structure, it is rather forbidding to study
obstructions in the category {\bf Graphs} directly. Instead, we
consider a~functor $\cf:\graphs\ra\ct$, where $\ct$ is some category
with a~well-developed obstruction theory.

For two graphs $A,B\in\co(\graphs)$, if there exists a~graph
homomorphism $\varphi:A\ra B$, then, since $\cf$ is a~functor, we have
an~induced morphism $\cf(\varphi):\cf(A)\ra\cf(B)$. If, on the other hand, by some general obstruction arguments in $\ct$, there can be no morphism
$\cf(A)\ra\cf(B)$, then we have gotten a~contradiction, hence
$\cm(A,B)=\emptyset$.

The question then becomes: how does one find ``good'' functors $\cf$,
i.e., functors which yield nontrivial obstructions to the existence of
morphisms in \graphs.

The centerpiece of this survey is the following choice of $\cf$:
choose a~test graph $T$, and map every $G\in\graphs$ to a~topological
space which is derived from the set of graph homomorphisms
$\varphi:T\ra G$. The idea of topologizing the set of graph
homomorphisms between two given graphs is due to L.\ Lov\'asz and will
be presented in detail in Section~\ref{sect3}.

Let us recall the following standard construction in category theory,
cf.\ \cite[Section II.6]{McL}. For a~category $\cc$ and an~object
$a\in\co(\cc)$, the {\em category of objects under $a$}, denoted
$a\downarrow\cc$, is defined as follows:

\nin $\bullet$ the objects of $a\downarrow\cc$ are all pairs $(m,b)$, where
$m$ is a~morphism from $a$ to $b$;

\nin $\bullet$ for $b_1,b_2\in\co(\cc)$, and $m_1\in\cm(a,b_1)$, $m_2\in\cm(a,b_2)$,
the morphisms in $a\downarrow\cc$ from $(m_1,b_1)$ to $(m_2,b_2)$ are
all morphisms $m:b_1\ra b_2$, such that $m\circ m_1=m_2$ (in other words, 
all ways to complement $m_1,m_2$ to an~appropriate commutative triangle).

The interesting and crucial detail here is the additional topological
structure which we have on top of the more usual comma category
construction.

\lecture{The functor $\thom(-,-)$.}

\section{Complexes of graph homomorphisms.} \label{sect3}

\subsection{Complex of complete bipartite subgraphs and 
the neighborhood complex.} \label{ss3.1} 


\nin First, we state the well-known definition in the generality which we need here.

\begin{df} \label{df:combip}
Let $A,B\subseteq V(G)$, $A,B\neq\emptyset$. We call $(A,B)$ a~{\bf
complete bipartite subgraph} of $G$, if for any $x\in A$, $y\in B$, we
have $(x,y)\in E(G)$, i.e., $A\times B\subseteq E(G)$.
\end{df}

In particular, note that all vertices in $A\cap B$ are required to have loops, and that the edges between the vertices of $A$ (or of $B$) are allowed.

\begin{figure}[hbt]
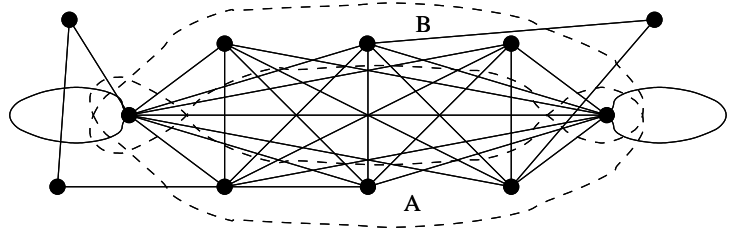

\begin{center}
  \begin{picture}(0,0)%
    \includegraphics{bip.pstex}%
  \end{picture}%
  \input{bip.pstex_t}%
 
\end{center}
\caption{A complete bipartite subgraph.}
\label{fig:bip}
\end{figure}

Let $G$ be a (possibly infinite) graph. Let $\Delta^{V(G)}$ be
a~simplex whose set of vertices is $V(G)$, in particular, the
simplices of $\Delta^{V(G)}$ can be identified with the finite subsets
of $V(G)$. We stress here, that we take as an~infinite simplex, the
colimit of the standard inclusion sequence of finite simplices:
\[\da^0\hra\da^1\hra\da^2\hra\dots.\] Under this convention, the points
of $\Delta^{V(G)}$ are all convex combinations of the points $V(G)$,
where only finitely many points have nonzero coefficients.

A~direct product of regular CW complexes is again a~regular CW
complex. Even stronger, $\Delta^{V(G)}\times\Delta^{V(G)}$ can be
thought of as a polyhedral complex, whose cells are direct products of
two simplices.

\begin{df} \label{df:bip}
  $\bip(G)$ is the subcomplex of $\Delta^{V(G)}\times\Delta^{V(G)}$
  defined by the following condition: $\sigma\times\tau\in\bip(G)$ if
  and only if $(\sigma,\tau)$ is a~complete bipartite subgraph of~$G$.
\end{df}

Note that if $(A,B)$ is a complete bipartite subgraph of $G$, and
$\wti A\subseteq A$, $\wti B\subseteq B$, $\wti A,\wti B\neq
\emptyset$, then $(\wti A,\wti B)$ is also a~complete bipartite
subgraph of~$G$. This verifies that $\bip(G)$ is actually
a~subcomplex.

$\bip(G)$ is a CW complex, whose closed cells are isomorphic to direct products of simplices (in the particular case here, they are in fact products of two simplices). We call complexes satisfying that property {\it prodsimplicial}.

In 1978 Lov\'asz proposed the following construction.

\begin{df} \label{df:nbd}
Let $G$ be a graph. The {\bf neighborhood complex} of $G$ is the
simplicial complex $\cn(G)$ defined as follows: its vertices are all
non-isolated vertices of~$G$, and its simplices are all the subsets of
$V(G)$ which have a~common neighbor.
\end{df}

Let $\cn(v)$ denote the set of neighbors of $v$, i.e.,
\[\cn(v)=\{x\in V(G)\,|\,(v,x)\in E(G)\}.\]
Then, the maximal simplices of $\cn(G)$ are precisely $\cn(v)$, for
$v\in V(G)$. The complexes $\bip(G)$ and $\cn(G)$ are closely related.

\begin{prop} \label{propn}
Let $G$ be an arbitrary graph.

\vskip3pt

\nin {\rm (a) (\cite[Proposition 4.2]{BK03b}).} 

\vskip3pt

\nin $\bip(G)$ is homotopy equivalent to~$\cn(G)$.

\vskip3pt

\nin {\rm (b) (\cite[Theorem 7.2]{K5}).} 

\vskip3pt

\nin Even stronger, $\bip(G)$ and $\cn(G)$ have the same 
simple homotopy type.
\end{prop}

$\bip(G)$ is our first example of the $\thom(-,-)$-construction,
namely, $\bip(G)$ is isomorphic, as a~polyhedral complex, to
$\thom(K_2,G)$.

\begin{figure}[hbt]
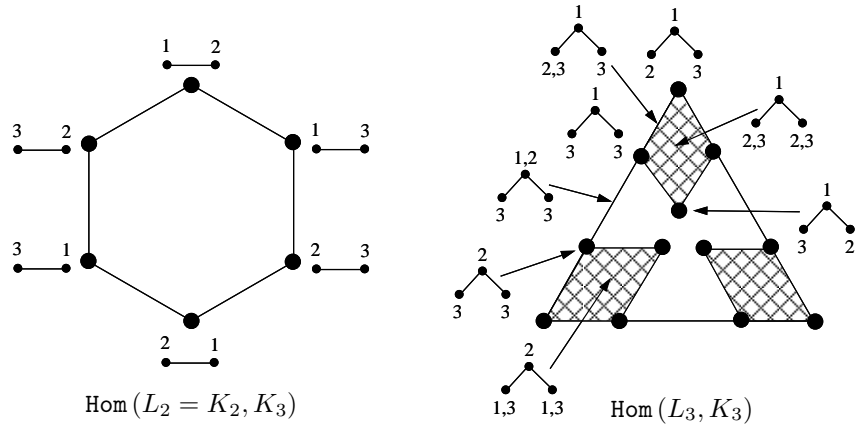

\begin{center}
  \begin{picture}(0,0)%
    \includegraphics{homlk2.pstex}%
  \end{picture}%
  \input{homlk2.pstex_t}%
 
\end{center}
\caption{3-coloring complexes of an edge and of a 3-string.}
\label{fig:homlk2}
\end{figure}

\subsection{$\thom$-construction for graphs.} \label{ss3.3} 


\nin We shall now define $\thom(T,G)$ for an arbitrary pair of graphs
$T$ and $G$. As a~model, we take the definition of $\bip(G)$. Let
again $\Delta^{V(G)}$ be a~simplex whose set of vertices is
$V(G)$. Let $C(T,G)$ denote the {\em weak } direct product
$\prod_{x\in V(T)} \Delta^{V(G)}$, i.e., the copies of $\Delta^{V(G)}$
are indexed by vertices of~$T$.

By the weak direct product we mean the following construction: a~cell
in $C(T,G)$ is a~direct product of cells $\prod_{x\in V(T)}\sigma_x$,
with the extra condition that $\dim\sigma_x\geq 1$ for only finitely
many~$x$. The dimension of this cell is $\sum_{x\in
V(T)}\dim\sigma_x$, in particular, it is finite.

\begin{df} \label{df:hom}
$\thom(T,G)$ is the subcomplex of $C(T,G)$ defined by the following
condition: $\sigma=\prod_{x\in V(T)}\sigma_x\in\thom(T,G)$ if and only
if for any $x,y\in V(T)$, if $(x,y)\in E(T)$, then
$(\sigma_x,\sigma_y)$ is a~complete bipartite subgraph of~$G$.
\end{df}

Let us make a number of simple, but fundamental, observations about the
complexes $\thom(T,G)$.

\nin (1) The topology of $\thom(T,G)$ is inherited from the product
topology of $C(T,G)$. By this inheritance, the cells of $\thom(T,G)$
are products of simplices.

\nin (2) $\thom(T,G)$ is a polyhedral complex whose cells are indexed
by all functions $\eta:V(T)\rightarrow
2^{V(G)}\setminus\{\emptyset\}$, such that if $(x,y)\in E(T)$, then
$\eta(x)\times\eta(y)\subseteq E(G)$.

The closure of a~cell $\eta$ consists of all cells indexed by
$\ti\eta:V(T)\rightarrow 2^{V(G)}\setminus\{\emptyset\}$, which
satisfy $\ti\eta(v)\subseteq\eta(v)$, for all $v\in V(T)$.

Throughout this survey, we shall make extensive use of the
$\eta$-notation.

\begin{figure}[hbt]
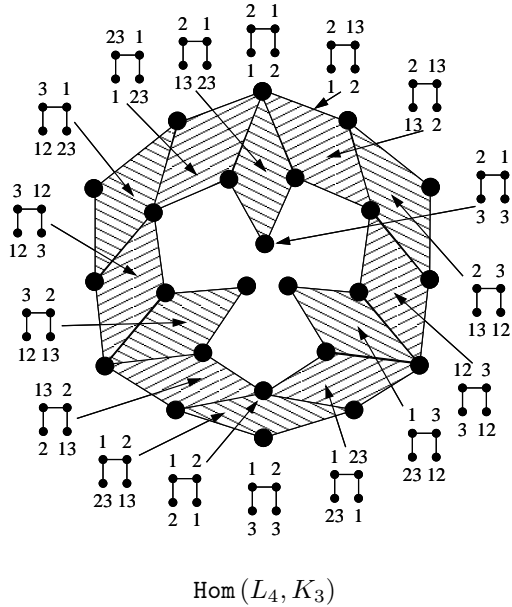

\begin{center}
  \begin{picture}(0,0)%
    \includegraphics{homlk3.pstex}%
  \end{picture}%
  \input{homlk3.pstex_t}%
 
\end{center}
\caption{3-coloring complex of a 4-string.}
\label{fig:homlk3}
\end{figure}

\nin (3) In the literature there are several different notations for
the set of all graph homomorphisms from a graph $T$ to the graph
$G$. Since an untangling of the definitions shows that this set is
precisely the set of vertices of $\thom(T,G)$, i.e., its 0-skeleton,
it feels natural to denote it by $\chom(T,G)$.

\nin (4) On the intuitive level, one can think of each $\eta:V(T)\ra
2^{V(G)}\setminus\{\emptyset\}$, satisfying the conditions of the
Definition~\ref{df:hom}, as associating non-empty lists of vertices of
$G$ to vertices of $T$ with the condition on this collection of lists
being that any choice of one vertex from each list will yield a~graph
homomorphism from $T$ to~$G$.

\nin (5) The standard way to turn a~polyhedral complex into
a~simplicial one is to take the barycentric subdivision. This is
readily done by taking the face poset and then taking its nerve (order
complex). So, here, if we consider the partially ordered set
$\cf(\thom(T,G))$ of all $\eta$ as in Definition~\ref{df:hom}, with
the partial order defined by $\ti\eta\leq\eta$ if and only if
$\ti\eta(v)\subseteq\eta(v)$, for all $v\in V(T)$, then we get that
the order complex $\Delta(\cf(\thom(T,G)))$ is a~barycentric
subdivision of $\thom(T,G)$. A~cell $\tau$ of $\thom(T,G)$ corresponds
to the union of all the simplices of $\Delta(\cf(\thom(T,G)))$ labeled
by the chains with the maximal element~$\tau$.

\vskip5pt

\begin{figure}[hbt]
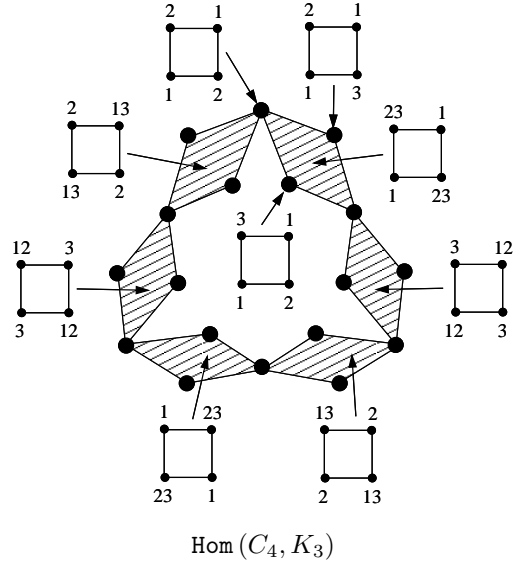

\begin{center}
  \begin{picture}(0,0)%
    \includegraphics{homc45k3.pstex}%
  \end{picture}%
  \input{homc45k3.pstex_t}%
 
\end{center}
\caption{3-coloring complex of a 4-cycle.}
\label{fig:c45k3}
\end{figure}

\nin Some examples are shown on Figures~\ref{fig:homlk2}, \ref{fig:homlk3}, \ref{fig:c45k3}, and \ref{fig:c5k3}. On these figures we used the following notations: $L_n$ denotes an
$n$-string, i.e., a~tree with $n$ vertices and no branching points,
$C_m$ denotes a~cycle with $m$ vertices, i.e., $V(C_m)=\dz_m$,
$E(C_m)=
\{(x,x+1),(x+1,x)\,| \,x\in\dz_m\}$.

\begin{figure}[hbt]
\begin{center}
  \begin{picture}(0,0)%
    \includegraphics{homc5k3.pstex}%
  \end{picture}%
  \input{homc5k3.pstex_t}%
 
\end{center}
\caption{3-coloring complex of a 5-cycle.}
\label{fig:c5k3}
\end{figure}

\section{Morphism complexes.} \label{sect3b}

\subsection{General construction.} \label{ss3.4} 


\nin As mentioned above, one way to interpret Definition~\ref{df:hom} is the following: the cells are indexed by the maps $\eta:V(T)\ra 2^{V(G)} \setminus\{\emptyset\}$, such that any choice $\varphi:V(T)\ra V(G)$, satisfying $\varphi(x)\in\eta(x)$, for all $x\in V(T)$, defines a~graph homomorphism $\varphi\in\thom_0(T,G)$. One can generalize this as follows.

Let $A$ and $B$ be two sets, and let $M$ be a~collection of some set
maps $\varphi:A\ra B$. Let $C(A,B)=\prod_{x\in A}\Delta^B$, where
$\Delta^B$ is the simplex having $B$ as a~vertex set, and copies in
the direct product are indexed by the elements of~$A$ (the direct
product is taken in the same weak sense as in the
subsection~\ref{ss3.3}).

\begin{df} \label{df:genhom}
Let $\thom_M(A,B)$ be the subcomplex of $C(A,B)$ consisting of all
$\sigma=\prod_{x\in A}\sigma_x$, such that any choice $\varphi:A\ra B$
satisfying $\varphi(x)\in\sigma_x$, for all $x\in A$, yields a~map in~$M$.
\end{df}

Intuitively one can think of the map $\varphi$ as the {\em section} of
$\sigma$, and the condition can then be verbally stated: {\em all
sections lie in~$M$.} Clearly, the complex $\thom_M(A,B)$ is always prodsimplicial.

An important fact is that $\thom_{-}(-,-)$ complexes are fully
determined by the low-dimensional data. This result did not previously
appear in the literature, so we include here a complete proof.
\begin{prop} \label{prop:1sk}
All complexes $\thom_{-}(-,-)$ with isomorphic 1-skeletons are
isomorphic to each other as polyhedral complexes. 

More precisely, the 1-skeleton determines the complex in the following
way: every product of simplices, whose 1-skeleton is in the 1-skeleton
of the complex, itself belongs to the complex.
\end{prop}

\pr Let us consider a complex $X$ of the type $\thom_{M}(A,B)$, 
where $A,B$ and $M$ are a priori unknown. Trivially, the 0-skeleton of
$X$ is the set $M$ itself. Furthermore, the 1-skeleton tells us which
pairs of set maps $\varphi,\psi:A\ra B$ differ precisely in one
element of~$A$.

Clearly, for $\sigma=\prod_{x\in A}\sigma_x\in C(A,B)$ to belong to
$\thom_M(A,B)$, it is required that the 1-skeleton of $\sigma$ is
a~subgraph of the 1-skeleton of $\thom_M(A,B)$. Let us show that the
converse of this statement is true as well.

Let $\Gamma$ be the 1-skeleton of $X$. For every edge $e$ in $\Gamma$
let $\lambda(e)\in A$ denote the element in which the value of the
function is changed along~$e$. Since we do not know the set $A$, we
can only make the statements of the type {\it the labels of these two
edges are the same/different}.

Assume we have $S\subseteq M$, such that $S$ can be written as
a~direct product $S=S_1\times S_2\times\dots\times S_t$. Assume
furthermore that the subgraph of $\Gamma$ induced by $S$ is precisely
the 1-skeleton of the corresponding cell.

First, consider 3 elements $a,b,c\in S_1\times\dots\times S_t$, which
have the same indices in all $S_i$'s except for exactly one. Then, by
our assumption on $S$, the subgraph of $\Gamma$ induced on the
vertices $a$, $b$, and $c$, is a~triangle. Clearly, if 3 changes of
a~value of a~function result in the same function, then the changes
were done in the same element of $A$, i.e.,
$\lambda(a,b)=\lambda(a,c)=\lambda(b,c)$.

Next, consider 4 elements $a,b,c,d\in S_1\times\dots\times S_t$, such
that pairs of vertices $(a,b),(b,c),(c,d)$, and $(a,d)$, have the same
indices in all $S_i$'s except for exactly one. Assume further that
this index is not the same for $(a,b)$ and $(b,c)$: say $a$ and $b$
differ in $S_1$, and $b$ and $c$ differ in $S_2$.
 
According to our assumption on $S$, $(a,b),(b,c),(c,d)$, and $(a,d)$
are edges of~$\Gamma$. If $\lambda(a,b)=\lambda(b,c)$, then $\Gamma$
contains the edge $(a,c)$, and $\lambda(a,b)=\lambda(a,c)$, which
contradicts our choice of~$S$.
 
If, on the other hand, $\lambda(a,b)\neq\lambda(b,c)$, then, since
changes of functions along the paths $a\ra b\ra c$ and $a\ra d\ra c$
should give the same answer, we are left with the only possibility:
namely $\lambda(a,b)=\lambda(c,d)$, and $\lambda(b,c)=\lambda(a,d)$.

Let $a\in S_1\times\dots\times S_t$, $a=(a_1,\dots,a_t)$. By our first
argument, if $b\in S_1\times\dots\times S_t$, $b=(a_1,\dots,\ti a_i,
\dots,a_t)$, then $\lambda(a,b)$ does not depend on $a_i$ and~$\ti
a_i$.
Furthermore, let $c,d\in S_1\times\dots\times S_t$,
$d=(a_1,\dots,\ti a_j,\dots,a_t)$, $c=(a_1,\dots,\ti a_i,\dots,\ti a_j,
\dots,a_t)$, for $i\neq j$. By our second argument, applied to
$a,b,c,d$, we get that $\lambda(a,b)=\lambda(c,d)$. If iterated for
various $j$, this implies that $\lambda(a,b)$ does not depend on
$a_1,\dots,a_{i-1},a_{i+1},\dots,a_t$ either; thus it may depend only
on the index~$i$.

Finally, this label should be different for different $i$'s, as
otherwise, by the same argument as above, we would get more edges in
the subgraph of $\Gamma$ induced by $S$, than what we allowed by our
assumptions.

Summarizing, we have shown that the cell $\sigma=\prod_{x\in A}\sigma_x\in
C(A,B)$ belongs to $\thom_M(A,B)$ if and only if the 1-skeleton of
$\sigma$ is a~subgraph of~$\Gamma$. This implies that $\thom_M(A,B)$
is uniquely determined by its 1-skeleton.  \qed

\vskip5pt

Intuitively, one can interpret Proposition~\ref{prop:1sk} as
saying that, with respect to its 1-skeleton, $\thom_M(A,B)$ is the
polyhedral analog of the flag complex construction.

\subsection{Specifying the parameters in the general construction.} \label{ss3.5} 


\nin (1) As mentioned above, if we take $A$ and $B$ to be the sets of
vertices of two graphs $T$ and $G$, and then take $M$ to be the set of
graph homomorphisms from $T$ to $G$, then $\thom_M(A,B)$ will coincide
with $\thom(T,G)$.

\vskip5pt

\nin (2) We think of a {\em directed} graph $G$ as a pair of sets
$(V(G),E(G))$, such that $E(G)\subseteq V(G)\times V(G)$.

\begin{df} \label{df:dirgrhom}
For two directed graphs $T$ and $G$, a {\bf directed graph
homomorphism} from $T$ to $G$ is a~map $\varphi:V(T)\ra V(G)$, such
that $(\varphi\times\varphi)(E(T))\subseteq E(G)$.
\end{df}

Let $A$ and $B$ be the sets of vertices of two directed graphs $T$ and
$G$, and let $M$ to be the set of directed graph homomorphisms from
$T$ to $G$, then $\thom_M(A,B)$ is the analog of $\thom(T,G)$ for
directed graphs.

For a~directed graph $G$, let $u(G)$ be the undirected graph obtained
from $G$ by forgetting the directions, and identifying the multiple
edges. We remark that for any two directed graphs $G$ and $H$, the
complexes $\thom(G,H)$ and $\thom(u(G),u(H))$ are isomorphic, if
$E(H)$ is $\zz$-invariant.

\begin{figure}[hbt]
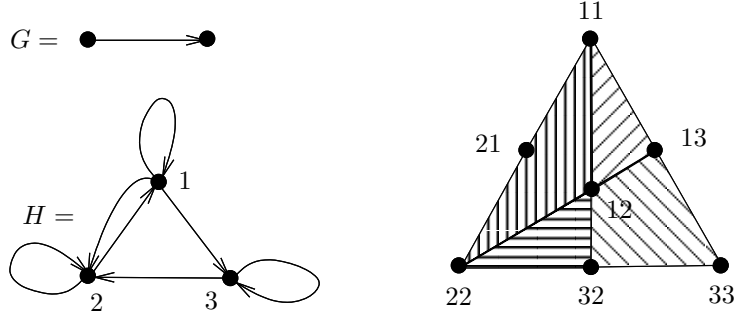

\begin{center}
  \begin{picture}(0,0)%
    \includegraphics{homdir.pstex}%
  \end{picture}%
  \input{homdir.pstex_t}%
 
\end{center}
\caption{An example of a $\thom$ complex for two directed graphs.}
\label{fig:homdir}
\end{figure}

\vskip5pt
 
\nin (3) Let $A$ and $B$ be the vertex sets of simplicial complexes
$\Delta_1$ and $\Delta_2$, and let $M$ be the set of simplicial maps
from $\Delta_1$ to $\Delta_2$, then $\thom_M(A,B)$ is the analog of
$\thom(T,G)$ for simplicial complexes.
 
\vskip5pt

\nin (4) Recall that a hypergraph with the vertex set $V$ is a subset
$\ch\subseteq 2^V$. Let $A$ and $B$ be the vertex sets of hypergraphs
$\ch_1$ and $\ch_2$. There are various choices for when to call
a~map $\varphi:A\ra B$ a~{\em hypergraph homomorphism}. Two
possibilities which we mention here are: one could require that
$\varphi(\ch_1)\subseteq\ch_2$, or one could ask that for any
$H_1\in\ch_1$, there exists $H_2\in\ch_2$, such that
$\varphi(H_1)\subseteq H_2$. The example (3) is a~special case of
both. Either way, the corresponding complex $\thom_M(A,B)$ provides us
with an~analog of $\thom(T,G)$ for hypergraphs.

\vskip5pt
 
\nin (5) Let $A$ and $B$ be the vertex sets of posets $P$ and $Q$, 
and let $M$ be the set of order-preserving maps from $P$ to $Q$, then
$\thom_M(A,B)$ is the analog of $\thom(T,G)$ for posets.
 
\section{Historic detour.} \label{sect4}

\subsection{The Kneser-Lov\'asz theorem.}

The Kneser conjecture was posed in 1955, see \cite{Knes}, and
concerned chromatic numbers of a specific family of graphs, later
called {\it Kneser graphs}. We denoted these graphs $K_{n,k}$, see
Definition~\ref{df:kngr}.  In 1978 L.\ Lov\'asz solved the Kneser
conjecture by finding geometric obstructions of Borsuk-Ulam type to
the existence of graph colorings.

\begin{thm} \label{thm:knlothm}
  {\rm (Kneser-Lov\'asz, \cite{Knes,Lo}).} 
  
  \nin For arbitrary positive integers $n,k$, such that $n\geq 2k$, we
  have $\chi(K_{n,k})=n-2k+2$.
\end{thm}

To show the inequality $\chi(K_{n,k})\geq n-2k+2$ Lov\'asz associated
the neighborhood complex $\cn(G)$, see Definition~\ref{df:nbd}, to
an~arbitrary graph $G$, and then used the connectivity information of
the topological space $\cn(G)$ to find obstructions to the
colorability of~$G$. More precisely, he proved the following
statement.

\begin{thm} \label{thm:lothm}
  {\rm (Lov\'asz, \cite{Lo}).}  

\nin {\it Let $G$ be a graph, such that
    $\cn(G)$ is $k$-connected for some $k\in\dz$, $k\geq -1$, then
    $\chi(G)\geq k+3$.}
\end{thm}

The main topological tool which Lov\'asz employed was the Borsuk-Ulam
theorem. Shorter proofs were obtained by B\'ar\'any, \cite{Bar}, and Greene, \cite{Gr}, both using some versions of the Borsuk-Ulam theorem, see also \cite{GoR}. A~nice brief survey of these can be found in~\cite{dL}.

Since these developments, the topological equivariant methods have
gained ground and became a~part of the standard repertoire in
combinatorics. We refer the reader to the series of papers
\cite{Ziv2,Ziv1,Ziv3} for an~excellent introduction to the subject.

We have seen in Proposition~\ref{propn}, that the complexes $\cn(G)$ and $\thom(K_2,G)$
have the same simple homotopy type. This fact leads one to consider the family of $\thom$ complexes as a natural context in which to look for further obstructions to the existence of graph homomorphisms. 

\subsection{Later developments.}

\subsubsection{The vertex-critical subgraphs of Kneser graphs.} \label{sss2.2.1} $\,$

\vskip5pt

\nin
For a graph $G$ and a vertex $v\in V(G)$ we introduce the following notation: $G-v$ denotes the graph which is obtained from $G$ by deleting the vertex $v$ and all edges adjacent to $v$, i.e., $V(G-v)=V(G)\sm\{v\}$, and $E(G-v)=E(G)\cap(V(G-v)\times V(G-v))$.

Shortly after Lov\'asz' result, Schrijver, in~\cite{Sch}, has sharpened Theorem~\ref{thm:knlothm}. To formulate his result, we recall that a~graph $G$ is called {\it vertex-critical} if, for any vertex $v\in V(G)$, we have $\chi(G)=\chi(G-v)+1$.

\begin{df}\label{df:stkngr}
Let $n,k$ be positive integers, $n\geq 2k$. The {\bf stable Kneser graph}
$K_{n,k}^\stab$ is defined to be the graph whose set of vertices is the set
of all $k$-subsets $S$ of $[n]$, such that if $i\in S$, then $i+1\notin S$, and, if $n\in S$, then $1\notin S$. Two subsets are joined by an edge if and only if they are disjoint.
\end{df}

Clearly, $K_{n,k}^\stab$ is an~induced subgraph of $K_{n,k}$. 

\begin{thm} \label{thm:schthm} {\rm (Schrijver, \cite{Sch}).}

\nin $K_{n,k}^\stab$ is a~vertex-critical subgraph of $K_{n,k}$, i.e., $K_{n,k}^\stab$ is a~vertex-critical graph, and $\chi(K_{n,k}^\stab)=n-2k+2$.
\end{thm}

\subsubsection{Chromatic numbers of Kneser hypergraphs.} \label{sss2.2.2} $\,$

\vskip5pt

\nin In 1986, Alon-Frankl-Lov\'asz, \cite{AFL}, have generalized Theorem~\ref{thm:knlothm} to the case of hypergraphs. To start with, recall the standard way to extend the notion of the chromatic number to hypergraphs. 

\begin{df} \label{df:chhyper}
For a~hypergraph $\ch$, the chromatic number $\chi(\ch)$ is, by definition, the minimal number of colors needed to color the vertices of $\ch$ so that no hyperedge is monochromatic. 
\end{df}

Next, there is a~standard way to generalize Definition~\ref{df:kngr} to the case of hypergraphs. 

\begin{df} \label{df:knhyp}
Let $n,k,r$ be positive integers, such that $r\geq 2$, and $n\geq rk$. The {\bf Kneser $r$-hypergraph} $K_{n,k}^r$ is the $r$-uniform hypergraph, whose ground set consists of all $k$-subsets of $[n]$, and the set of hyperedges consists of all $r$-tuples of disjoint $k$-subsets.
\end{df}

Using the introduced notations, we can now formulate the generalization of Theorem~\ref{thm:knlothm}. 

\begin{thm} \label{thm:hyperkn} {\rm (Alon-Frankl-Lov\'asz, \cite{AFL}).}

\nin 
  \nin For arbitrary positive integers $n,k,r$, such that $r\geq 2$, and $n\geq rk$, we have
\[\chi(K_{n,k}^r)=\left\lceil \frac{n-rk+r}{r-1}\right\rceil.
\]
\end{thm}

Theorem~\ref{thm:hyperkn} can be proved using the generalization of the Borsuk-Ulam theorem from \cite{BSS}.

\subsubsection{Further references.} \label{sss2.2.3} $\,$

\vskip5pt

\nin There has been a~substantial body of further important work, which, due to space constraints, we do not pursue in detail in this survey, some of the references are \cite{Do,Kr,Kr2,Ma2,MZ,Sa,Zi02}.

There have also been multiple constructions, such as box complexes,
designed to generalize the original Lov\'asz neighbourhood
complexes. However, as later research showed, the bounds obtained in
that way were essentially convertible, since the $\zz$-homotopy types
of these complexes were very closely related, either by simply being
the same, or by means of one being the suspension of another, or
something close to that. This means that all these constructions are
avatars of the same object, as explained in~\cite{Ziv4}.

\section{More about the $\thom$-complexes. }

\subsection{Coproducts.} \label{sss5.1.1}


\nin For any three graphs $G$, $H$, and $K$, we have
\begin{equation} \label{eq:coprod}
 \thom(G\coprod H,K)=\thom(G,K)\times\thom(H,K),
\end{equation}
and, if $G$ is connected, and $G\neq K_1$, then also
\[\thom(G,H\coprod K)=\thom(G,H)\coprod\thom(G,K),\]
where the equality denotes isomorphism of polyhedral complexes.

The first formula is obvious. To verify the second one, note that, for any graph homomorphism $\eta:V(G)\ra 2^{V(H)\cup V(K)}\sm\{\emptyset\}$, and any $x,y\in V(G)$, such that $(x,y)\in E(G)$, if $\eta(x)\cap V(H)\neq \emptyset$, then $\eta(y)\subseteq V(H)$, which under the assumptions on $G$ implies that $\bigcup_{x\in V(G)}\eta(x)\subseteq V(H)$.

\subsection{Products.} \label{sss5.1.1b}

\nin For any three graphs $G$, $H$, and $K$, we have the following homotopy equivalence, see \cite{Ba05}:
\begin{equation} \label{eq:prod1}
 \thom(G,H\times K)\simeq\thom(G,H)\times\thom(G,K).
\end{equation}

In fact, the formula \eqref{eq:prod1} can be strengthened to state that the left hand side is simple homotopy equivalent (in the sense of Whitehead, see \cite{Co73}) to the right hand side. Since this simple homotopy equivalence result is new, we include a~complete argument, as promised in the abstract. 

Consider the following three maps $2^{p_H}:2^{V(H)\times V(K)}\ra 2^{V(H)}$, $2^{p_K}:2^{V(H)\times V(K)}\ra 2^{V(K)}$, and $c:2^{V(H)}\times 2^{V(K)}\ra 2^{V(H)\times V(K)}$, where $2^{p_H}$ and $2^{p_K}$ are induced by the standard projection maps $p_H:V(H)\times V(K)\ra V(H)$ and $p_K:V(H)\times V(K)\ra V(K)$, and $c$ is given by $c(A,B)=A\times B$.
We let $\psi:2^{V(H)\times V(K)}\ra 2^{V(H)\times V(K)}$ denote the composition map $\psi(S)=c(2^{p_H}(S),2^{p_K}(S))=2^{p_H}(S)\times 2^{p_K}(S)$.

Given a~cell of $\thom(G,H\times K)$ indexed by $\eta:V(G)\ra 2^{V(H)\times V(K)}\sm\{\emptyset\}$, one can see that the composition function $\psi\circ\eta:V(G)\ra 2^{V(H)\times V(K)}\sm\{\emptyset\}$ will also index a~cell. Indeed, for any $(x,y)\in E(G)$ we know that $(\eta(x),\eta(y))$ is a~complete bipartite subgraph of $H\times K$, which is the same as to say that, for any $\alpha\in\eta(x)$, and $\beta\in\eta(y)$, we have $(p_H(\alpha),p_H(\beta))\in E(H)$, and $(p_K(\alpha),p_K(\beta))\in E(K)$. If we now choose $\ti\alpha\in\psi(\eta(x))$, and $\ti\beta\in\psi(\eta(y))$, we have $p_H(\ti\alpha)=p_H(\alpha)$, for some $\alpha\in\eta(x)$, and $p_H(\ti\beta)=p_H(\beta)$, for some $\beta\in\eta(y)$, hence verifying that $(p_H(\ti\alpha),p_H(\ti\beta))\in E(H)$. The fact that $(p_K(\ti\alpha),p_K(\ti\beta))\in E(K)$ can be proved analogously.

This means that we have a~map $\varphi:\cf(\thom(G,H\times K))\ra \cf(\thom(G,H\times K))$. It is easy to see that $\varphi$ is order-preserving and ascending (meaning $\varphi(x)\geq x$, for any $x\in\cf(\thom(G,H\times K))$). It follows from \cite[Theorem 3.1]{K5a} that $\da(\cf(\thom(G,H\times K)))=\bd(\thom(G,H\times K))$ collapses onto $\da(\im\varphi)$. 

On the other hand, $\cf(\thom(G,H))\times\cf(\thom(G,K))\cong\im\varphi$ with the isomorphism given by the map $(\eta_1,\eta_2)\mapsto\eta$, where $\eta(x)=\eta_1(x)\times\eta_2(x)$, for any $x\in V(G)$. Thus we conclude that $\bd(\thom(G,H\times K))$ collapses onto $\da(\cf(\thom(G,H))\times\cf(\thom(G,K)))\cong\da(\cf(\thom(G,H)))\times\da(\cf(\thom(G,K)))=\bd(\thom(G,H))\times\bd(\thom(G,K))\cong\thom(G,H)\times\thom(G,K)$, and our argument is now complete.

\vskip5pt

For the analog of the formula \eqref{eq:prod1}, where the direct product is taken on the left, we need the following additional standard notion.

\begin{df} \label{df:grpow}
For two graphs $H$ and $K$, the {\bf power graph} $K^H$ is defined by
\begin{itemize}
	\item $V(K^H)$ is the set of all set maps $f:V(H)\ra V(K)$;
	\item $(f,g)\in E(K^H)$, for $f,g:V(H)\ra V(K)$, if and only if, whenever $(v,w)\in E(H)$, we also have $(f(v),g(w))\in E(K)$.
\end{itemize}
\end{df}

It is easy to see that the power graph notion is introduced precisely so that for any triple of graphs the following adjunction relation holds:
\begin{equation} \label{eq:adj}
 \chom(G\times H,K)=\chom(G,K^H).
\end{equation}

In our topological situation the formula \eqref{eq:adj} generalizes up to homotopy. More precisely, we have the following homotopy equivalence, see \cite{Ba05},
\begin{equation} \label{eq:prod2}
 \thom(G\times H,K)\simeq\thom(G,K^H).
\end{equation}

The formula \eqref{eq:prod2} can as well be strengthened to yield a~simple homotopy equivalence. Below we include a~complete argument. 

Define a~map $\psi:2^{V(K^H)}\ra 2^{V(K^H)}$, $\psi:\Omega\mapsto\psi(\Omega)$, as follows: $g\in\psi(\Omega)$ if and only if $g(x)\in\{f(x)\,|\,f\in\Omega\}$, for all $x\in V(H)$. In other words, we use the collection of functions $\Omega$ to specify the sets of values, which functions from $\psi(\Omega)$ are allowed to take. Clearly, we have $\psi(\Omega)\supseteq\Omega$. Take a~cell of $\thom(G,K^H)$, $\eta:V(G)\ra 2^{V(K^H)}\sm\{\emptyset\}$, and consider the composition map $\psi\circ\eta:V(G)\ra 2^{V(K^H)}\sm\{\emptyset\}$. Since $\eta$ is a~cell, we know that if $(x,y)\in E(G)$, and $\alpha\in\eta(x)$, $\beta\in\eta(y)$, then $(\alpha,\beta)\in E(K^H)$, i.e., whenever $(v,w)\in E(H)$, we have $(\alpha(v),\beta(w))\in E(K)$.

Choose $\ti\alpha\in\psi(\eta(x))$, and $\ti\beta\in\psi(\eta(y))$. To check that $(\ti\alpha,\ti\beta)\in E(K^H)$, we need to check that for any $(v,w)\in E(H)$, we have $(\ti\alpha(v),\ti\beta(w))\in E(K)$. However, by the definition of $\psi$, we know that $\ti\alpha(v)=\alpha(v)$, for some $\alpha\in\eta(x)$, and $\ti\beta(w)=\beta(w)$, for some $\beta\in\eta(y)$. It follows that $(\ti\alpha(v),\ti\beta(w))=(\alpha(v),\beta(w))\in E(K)$, and hence $\psi\circ\eta$ is again a~cell.

As a consequence, the composition gives us an order-preserving ascending map $\varphi:\cf(\thom(G,K^H))\ra\cf(\thom(G,K^H))$. The image of this map is isomorphic to $\cf(\thom(G\times H,K))$. The isomorphism map takes the poset element $\eta:V(G)\times V(H)\ra 2^{V(K)}\sm\{\emptyset\}$ to the poset element $\ti\eta:V(G)\ra 2^{V(K^H)}\sm\{\emptyset\}$ defined by
\[\ti\eta(x)=\{f:V(H)\ra V(K)\,|\,f(v)\in\eta(x,v),\text{ for all }v\in V(H)\},
\]
for all $x\in V(G)$. By \cite[Theorem 3.1]{K5a}, we conclude that the complex $\da(\cf(\thom(G,K^H)))=\bd(\thom(G,K^H))$ collapses onto its subcomplex $\da(\im\varphi)=\bd(\thom(G\times H,K))$.

We obtain an interesting special case of the formula \eqref{eq:prod2}
when substituting $G=K_1^o$ (which means a~graph with one looped
vertex).  Since $K_1^o\times H=H$, for any graph $H$, we conclude that
$\thom(H,K)\simeq\thom(K_1^o,K^H)$ for any two graphs $H$ and~$K$.  As
seen directly, for an arbitrary graph $G$, $\thom(K_1^o,G)$ is the
clique complex of the looped part of $G$, i.e., of the subgraph
induced by the set of vertices which have loops. In particular, the
complex $\thom(K_1^o,G)$ is simplicial. On the other hand, a~vertex
$f\in V(K^H)$ has a~loop if and only if $f$ is a~graph homomorphism.
We can therefore conclude that for arbitrary graphs $H$ and~$K$ the
complex $\thom(H,K)$ is homotopy equivalent to the clique complex of
the subgraph of $K^H$, induced by the set of all graph homomorphisms
from $H$ to~$K$.

\subsection{Associated covariant and contravariant functors.} \label{sss5.1.2} 
\label{ss:231}


\nin For an arbitrary polyhedral complex $X$, we let $\cf(X)$ denote its face
poset ordered by inclusion. The notion of a~{\it link} of a~vertex of
a~polyhedral complex allows several interpretations, so let us fix our
convention here. Let $v$ be a~vertex of $X$, the link of $v$, denoted
$\link_X(v)$, is the cell complex whose face poset is given by
$\cf_{>v}(X)$. Geometrically, $\link_X(v)$ can be obtained as
follows. Realize faces of $X$ as polyhedra, in a~coherent manner, and
take $\varepsilon$ to be a~positive number which is smaller then the
minimal length of an edge from $v$. Each face $F$ containing $v$ can
be truncated at vertex $v$, by cutting it along the set of points at
distance $\varepsilon$ from $v$. These cuts fit coherently to form the
desired cell complex. For example, the link of any vertex of a~cube is
a~triangle, and in general, the link of any vertex of a~polytope $K$
is the polytope obtained by truncating $K$ at the vertex~$v$.

Let $T,G$, and $K$ be three arbitrary graphs, and let $\varphi$ be
a~graph homomorphism from $G$ to $K$. Then the composition induces
a~poset map $f:\cf(\thom(T,G))\ra\cf(\thom(T,K))$, namely, for
$\eta:V(T)\ra 2^{V(G)}\sm\{\emptyset\}$, we have
$f(\eta)=2^\varphi\circ\eta$, where $2^\varphi$ is the map induced on
the subsets.

Recall that for arbitrary regular CW complexes $A$ and $B$, a poset
map $f:\cf(A)\ra\cf(B)$ comes from a~cellular map $\varphi:A\ra B$
(meaning that $f=\cf(\varphi)$), if and only if $\rk\varphi(x)\leq\rk
x$, for all $x\in\cf(A)$.  It is not difficult to check that this
condition is satisfied by the poset map
$f:\cf(\thom(T,G))\ra\cf(\thom(T,K))$ defined above. Hence, we can
conclude that this $f$ comes from a~cellular map from $\thom(T,G)$ to
$\thom(T,K)$, which we denote by~$\varphi^T$.

Moreover, a~detailed pointwise analysis of the polyhedral structure of $\thom(T,G)$ shows that cells (direct products of simplices) map surjectively to other cells, and that this map is a~product map induced by the corresponding maps on the simplices. Therefore, $\varphi^T$ is a~polyhedral map.

The situation is slightly more complicated if one considers the
functoriality in the first argument. Let us choose some {\it proper} 
graph homomorphism $\psi$ from $T$ to~$G$, and let $K$ be some graph. 
Again, by using composition we can define a poset map
$g:\cf(\thom(G,K))\ra\cf(\thom(T,K))$, namely, for $\eta:V(G)\ra
2^{V(K)}\sm\{\emptyset\}$, and $v\in V(T)$, we have
$g(\eta)(v)=\eta(\psi(v))$. This map is well-defined, since, first if
$v,w\in V(T)$, and $(v,w)\in E(T)$, then $(\psi(v),\psi(w))\in E(G)$,
and therefore, for any $x\in\eta(\psi(v))$, and $y\in\eta(\psi(w))$,
we have $(x,y)\in E(K)$, and second, by the properness assumption,
$\sum_{v\in V(T)}(|\eta(\psi(v))|-1)<\infty$. Furthermore, this map is
order-preserving: if $\tau\geq\eta$, i.e., if
$\tau(w)\supseteq\eta(w)$, for any $w\in V(T)$, then
$g(\tau)(w)=\tau(\psi(w))\supseteq\eta(\psi(w))=g(\eta)(w)$.

Intuitively, one can think of the map $g$ as the pullback map. It is important to remark that, if $\psi$ is not injective, it may happen that $\dim g(\eta) >\dim\eta$.

For an arbitrary regular CW complex $X$, let $\bd(X)$ denote the {\it
barycentric subdivision} of $X$. Since $g$ is an order-preserving map,
the induced map \[\Delta(g):\bd(\thom(G,K))\ra\bd(\thom(T,K))\] is
simplicial and gives the corresponding map of topological spaces,
which we denote~$\psi_K$.  However, $g$ does not always come from
a~cellular map. In fact, one can check that there exists a~cellular
map $\psi_K:\thom(G,K)\ra\thom(T,K)$, such that $\cf(\psi_K)=g$, if
$\psi$ is injective on the vertices of~$T$.

In any case, we see that $\thom(T,-)$ is a covariant functor from {\bf
Graphs} to {\bf Top}, while $\thom(-,K)$ is a~contravariant functor
from {\graphsb} to {\bf Top}; here {\bf Top} denotes the category
whose objects are topological spaces, and whose morphisms are all
continuous maps.

\subsection{Composition of $\,{\tt Hom}$'s.} \label{sss5.1.3} 


\nin For three arbitrary graphs $T,G$, and $K$, there is a~composition map 
\[\xi:\cf(\thom(T,G))\times\cf(\thom(G,K))\lra\cf(\thom(T,K)),\]
whose detailed description is as follows: for graph homomorphisms $\alpha:V(T)\ra 2^{V(G)}\sm\{\emptyset\}$, and $\beta:V(G)\ra 2^{V(K)}\sm\{\emptyset\}$, define the map $\ti\beta:2^{V(G)}\sm\{\emptyset\}\ra 2^{V(K)}\sm\{\emptyset\}$ by
\[\text{for } S\in 2^{V(G)}\sm\{\emptyset\}, \quad \ti\beta(S):=\cup_{x\in S}\beta(x),\qquad \quad
\]
and then set $\xi(\alpha,\beta):=(\ti\beta\circ\alpha:V(T)\ra
2^{V(K)}\sm\{\emptyset\})$. It is easy to check that this map is
well-defined. Indeed, let $x,y\in V(T)$, such that $(x,y)\in E(T)$,
choose arbitrary $a\in\alpha(x)$, and $b\in\alpha(y)$, and then choose
arbitrary $\ti a\in\beta(a)$, and $\ti b\in\beta(b)$. Clearly,
$(x,y)\in E(T)$ implies $(a,b)\in E(G)$, since $\alpha$ is a graph
homomorphism, which then implies $(\ti a,\ti b)\in E(K)$, since
$\beta$ is a graph homomorphism.

Applying the nerve functor $\Delta$ to the poset map $\xi$ we get a~simplicial map
\[\Delta(\xi):\da(\cf(\thom(T,G))\times\cf(\thom(G,K)))\lra\bd(\thom(T,K)),
\]
and hence, since for any posets $P_1$ and $P_2$, the simplicial complex $\da(P_1\times P_2)$ is homeomorphic to the polyhedral complex $\da(P_1)\times\da(P_2)$ (in fact it is its subdivision), we have a corresponding topological map
	\[\thom(T,G)\times\thom(G,K)\lra\thom(T,K).
\]
         
\subsection{Action of automorphism groups.} \label{sss5.1.3b} 


\nin An important special case of the situation described in the Subsection~\ref{sss5.1.3} is when the considered graph homomorphisms are actually isomorphisms. In other words, for arbitrary graphs $T$ and $G$, the elements $\varphi\in\aut(T)$ and $\psi\in\aut(G)$ induce polyhedral maps $\varphi_G:\thom(T,G)\ra\thom(T,G)$ and $\psi^T:\thom(T,G)\ra\thom(T,G)$, which are easily shown to be isomorphisms.

Summarizing, we have a~polyhedral action of the group $\aut(T)\times\aut(G)$ on the polyhedral complex $\thom(T,G)$. As an example, we have $\cs_m\times\cs_n$-action on $\thom(K_m,K_n)$, and $\cd_m\times\cs_n$-action on $\thom(C_m,K_n)$, where $\cs_n$ is the $n$-th symmetric group, and  $\cd_n$ is the $n$-th dihedral group.
         
We note the following useful fact: if for some vertex $v$ there exists
a~group element $\varphi\in\aut(T)$, such that $(v,\varphi(v))\in
E(T)$ (for example, if $\varphi$ flips an edge in $T$), then the
induced map $\varphi_G:\thom(T,G)\ra\thom(T,G)$ is fixed-point free
for an arbitrary graph $G$ without loops. For example, $\zz$
(corresponding to an arbitrary reflection from $\cd_{2r+1}$) acts
freely on $\thom(C_{2r+1},K_n)$.
         
\subsection{Universality.} \label{sss5.1.4} 


\nin In topological combinatorics it happens very often that the
family of the studied objects is universal with respect to the
invariants which one is interested in computing. This is also the case
not only for the $\thom$-complexes, but even for the
$\thom(K_2,-)$-complexes. The following result is due to Csorba,
\cite{Cs1}, and, independently, to \v{Z}ivaljevi\'{c}, \cite{Ziv4}.

\begin{thm} \label{thm:zivcs}
{\rm (\cite{Cs1,Ziv4}).}

\nin For each finite, free $\zz$-complex $X$, there exists a~graph $G$,
such that $\thom(K_2,G)$ is $\zz$-homotopy equivalent to~$X$.
\end{thm}

We note that Theorem~\ref{thm:zivcs} can be verified by combining
\cite[Theorem~32]{Ziv4}, with a~remark in the beginning of
\cite[Section~7]{Ziv4}.

\section{Folds.} \label{ss5.2} 

\subsection{Sequences of collapses induced by folds.}


\nin $\thom$-complexes behave well with respect to the following standard operation from graph
theory.

\begin{df} \label{df:fold}
For a graph $G$ and $v\in V(G)$, $G-v$ is called a {\bf fold} of $G$
if there exists $u\in V(G)$, $u\neq v$, such that $\cn(u)\supseteq\cn(v)$.
\end{df}

Let $G-v$ be a fold of $G$. We let $i:G-v\hookrightarrow G$ denote the
inclusion homomorphism, and let $f:G\ra G-v$ denote the folding
homomorphism defined by 
\[f(x)=\begin{cases} 
u,& \text{ for } x=v;\\
x,& \text{ for } x\neq v.
\end{cases}\]

Note that $i$ is a graph homomorphism for an arbitrary choice of $v\in V(G)$, whereas $f$ is a~graph homomorphism if and only if $G-v$ is a~fold,
in particular, this could be taken as an alternative definition of the fold. 

Let $X$ be a~polyhedral complex. Recall that an {\em elementary
collapse} of $X$ is a~removal of a pair of polyhedra $(\sigma,\tau)$,
such that $\sigma$ is a maximal polyhedron, $\dim\tau=\dim\sigma-1$,
and $\sigma$ is the only maximal polyhedron containing
$\tau$. Furthermore, let $Y$ be a~subcomplex of $X$. We say that $X$
{\it collapses onto} $Y$ if there exists a~sequence of elementary
collapses leading from $X$ to $Y$. If $X$ collapses onto $Y$, then $Y$
is a~strong deformation retract of~$X$.

\begin{figure}[hbt]
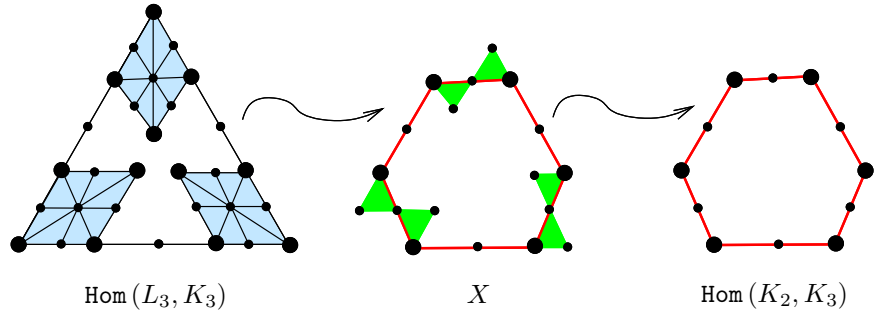

\begin{center}
  \begin{picture}(0,0)%
    \includegraphics{coll1.pstex}%
  \end{picture}%
  \input{coll1.pstex_t}%
 
\end{center}
\caption{A two-step folding of the first argument in $\thom(L_3,K_3)$.}
\label{fig:coll1}
\end{figure}

\begin{thm}\label{thm3.3k4} {\rm(\cite[Theorem 3.3]{K4}).} 

\nin Let $G-v$ be a fold of $G$ and let $H$ be some graph. Then
\begin{enumerate} 
 \item $\bd\thom(G,H)$ collapses onto $\bd\thom(G-v,H)$;
 \item $\thom(H,G)$ collapses onto $\thom(H,G-v)$. 
\end{enumerate}
The maps $i_H$ and $f^H$ are strong deformation retractions.
\end{thm}

Figure \ref{fig:coll1} shows an example of the collapsing sequence appearing in the proof of Theorem~\ref{thm3.3k4} (1).

\begin{figure}[hbt]
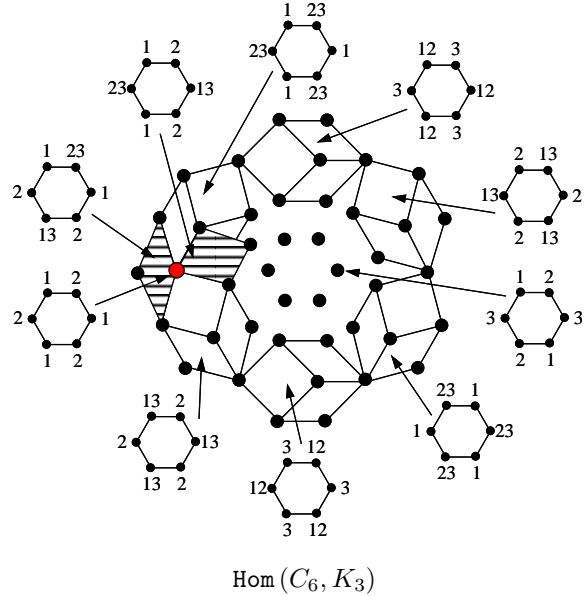

\begin{center}
  \begin{picture}(0,0)%
    \includegraphics{homc6k3.pstex}%
  \end{picture}%
  \input{homc6k3.pstex_t}%
 
\end{center}
\caption{The figure depicts the polyhedral complex of all 3-colorings of a~6-cycle. }
\label{fig:homc6k3}
\end{figure}

\begin{figure}[hbt]
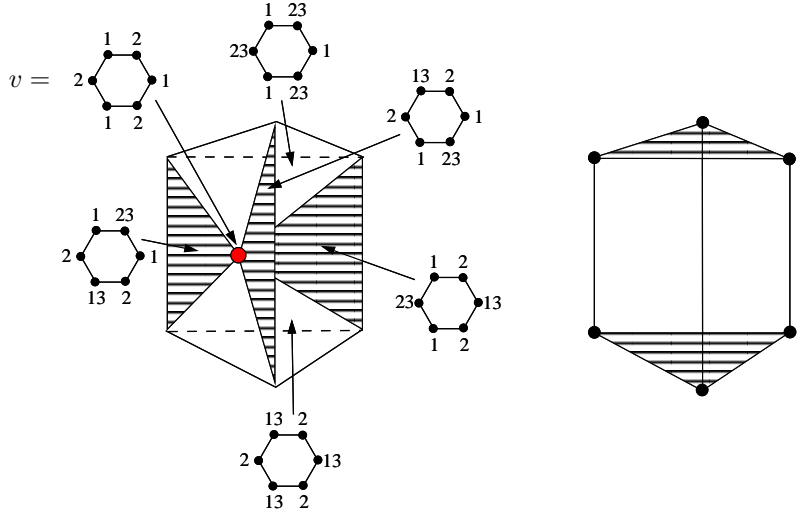

\begin{center}
  \begin{picture}(0,0)%
    \includegraphics{homc6k3b.pstex}%
  \end{picture}%
  \input{homc6k3b.pstex_t}%
 
\end{center}
\caption{
The figure on the left shows the neighbourhood of the vertex $v$ of
$\thom(C_6,K_3)$.  The figure on the right shows the link of this
vertex.}
\label{fig:homc6k3b}
\end{figure}

Note that Theorem~\ref{thm3.3k4} cannot be generalized to
encompass arbitrary graph homomorphisms $\phi$ of $G$ onto $H$, where
$H$ is a~subgraph of $G$, and $\phi$ is identity on $H$. For example,
$\thom(C_6,K_3)\not\simeq\thom(K_2,K_3)$, see
Figures~\ref{fig:homc6k3} and~\ref{fig:homc6k3b}, despite of the
existence of the sequence of graph homomorphisms $K_2\ra C_6\ra K_2$
which compose to give an~identity.

We remark that, for the sake of transparency, the striped rectangles
are shown on Figure~\ref{fig:homc6k3} only around one of the 6 joining
vertices, and only two out of the three. The big connected component
corresponds to the graph homomorphisms $\varphi:C_6\ra K_3$ having the
winding number~$0$. The isolated points correspond to the 6 possible
tight windings of $C_6$ around $K_3$. Observe also that the cubes are
solid.

\subsection{Applications.}

When $G$ is a graph, and $H$ is an induced subgraph of $G$, we say that $G$ {\it reduces to} $H$, if there exists a~sequence of folds leading from $G$ to~$H$.

\pagebreak

\begin{crl}\label{equicrl} {\rm(\cite[Corollary 5.3]{BK03b}).}

\nin Let $G$ be a~graph, and $S\subseteq V[G]$, such that $G$ reduces 
to $G[S]$. Assume $S$ is $\Gamma$-invariant for some
  $\Gamma\subseteq\aut(G)$. Then the inclusion $i:G[S]\hra G$ induces
  a~$\Gamma$-invariant homotopy equivalence
  $i_H:\thom(G,H)\ra\thom(G[S],H)$ for an~arbitrary graph~$H$.
\end{crl}

The Theorem \ref{thm3.3k4} can be used to obtain complete understanding of the homotopy type of the $\thom$-complexes for certain specific families of graphs.

\begin{prop} \label{pr:tree}
If $T$ is a~tree with at least one edge, and $G$ an arbitrary graph, then $\thom(T,G)$ is homotopy equivalent to $\cn(G)$. As a~consequence, if $F$ is a~forest, and $T_1,\dots,T_k$ are all its connected components consisting of at least 2 vertices, then $\thom(F,G)\simeq\prod_{i=1}^k\cn(G)$.
\end{prop}

An even more special case was important in \cite{BK03b,BK03c} for the proof of Lov\'asz Conjecture.

\begin{crl} \label{crl:tree}  {\rm(\cite[Proposition 5.4]{BK03b}).}

\nin If $T$ is a~finite tree with at least one edge, then the map $i_{K_n}:\thom(T,K_n)\ra\thom(K_2,K_n)$ induced by any inclusion $i:K_2\hookrightarrow T$ is a~homotopy equivalence, in particular $\thom(T,K_n)\simeq S^{n-2}$.
  
\nin If $F$ is a~finite forest, and $T_1,\dots,T_k$ are all its connected components consisting of at least 2 vertices, then $\thom(F,K_n)\simeq\prod_{i=1}^k S^{n-2}$.
\end{crl}

In this case, Corollary~\ref{equicrl} can be applied to describe the $\zz$-homotopy type as well. First, some new notations: let $S_a^n$, resp.\ $S_t^n$, denote the $n$-dimensional sphere, equipped with an~antipodal, resp.\ trivial $\zz$-action, where $n$ is a~nonnegative integer, or infinity.

\begin{prop} \label{pr:tree2}  {\rm(\cite[Proposition 5.5]{BK03b}).}

\nin Let $T$ be a~finite tree with at least one edge and a~$\zz$-action determined by an~invertible graph homomorphism $\gamma:T\ra T$. If $\gamma$ flips an~edge in $T$, then $\thom(T,K_n)\simeq_\zz S^{n-2}_a$, otherwise $\thom(T,K_n)\simeq_\zz S^{n-2}_t$.
\end{prop}

Curiously, another computable special case is that of an~unlooped complement of a~forest.

\begin{prop} \label{homcompf}
  Let $F$ be a~finite forest, and let $G$ be an~arbitrary graph, then $\thom(\ovr{F},G)\simeq\thom(K_m,G)$, where $m$ is the maximal cardinality of an~independent set in~$F$. 
\end{prop}

In particular, as was shown in \cite[Proposition 5.6]{BK03b} that
$\thom(\ovr{F},K_n)$ is homotopy equivalent to $\thom(K_m,K_n)$, and
hence, by Theorem~\ref{thm:kmkn}, to a wedge of
$(n-m)$-dimensional spheres.

\begin{rem}
Recently, folds gained further prominence in connection with $\thom$ complexes. One has discovered, see \cite{Ba05}, that it is possible to introduce a~natural Quillen model category structure, see \cite[Chapter V]{GeMa}, and \cite{Q}, on the category of graphs, such that the weak equivalences are precisely the maps, which allow factorizations into a~sequence of folds and unfolds (which therefore may be viewed as trivial homotopy equivalences).
\end{rem}

\lecture{Stiefel-Whitney classes and first applications.} 
 
\section{Elements of the principal bundle theory.}

\subsection{Spaces with a free action of a finite group and special cohomology elements.} \label{ssect:gammaspaces} 


\nin Consider a~regular CW complex $X$ with a~cellular action 
of a~finite group $\Gamma$. If desired, the $\Gamma$-action can be made to be simplicial by passing to the barycentric subdivision, (cf.\ \cite{Bre,Hat}).

For the interested reader we remark here that sometimes one takes the
barycentric subdivision even if the original action already was
simplicial. The main point of this is that one can make the action
enjoy an~additional property: {\em if a~simplex is preserved by one of
the group elements, then it must be pointwise fixed by this element.}

Next, assume that $\Gamma$ acts freely. In this case $X$ is called
a~{\em $\Gamma$-space}. We know that, by the general theory of
principal $\Gamma$-bundles, see e.g., \cite{tD}, there exists
a~$\Gamma$-equivariant map $w:X\ra{\bf E}\Gamma$, and, that the
induced quotient map $w/\Gamma:X/\Gamma\ra{\bf E}\Gamma/\Gamma={\bf
B}\Gamma$ is unique up to homotopy. Here ${\bf E}\Gamma$ denotes
a~contractible space on which the group $\Gamma$ acts freely, and
${\bf B}\Gamma$ denotes the {\it classifying space} of $\Gamma$ (also
known as the associated Eilenberg-MacLane space, or
$K(\Gamma,1)$-space, see e.g., \cite{AM,Bre93,GeMa,Hat,McL2,May,Wh78}).

Passing to cohomology, we see that the induced map
$(w/\Gamma)^*:H^*({\bf B}\Gamma)\ra H^*(X/\Gamma)$ does not depend on
the choice of $w$, and thus, the image of $(w/\Gamma)^*$ consists of
some canonically distinguished cohomology elements. For $z\in H^*({\bf
B}\Gamma)$, we let $w(z,X)$ denote the element $(w/\Gamma)^*(z)$,
which we call {\it characteristic class} associated to~$z$.

\begin{figure}[hbt]
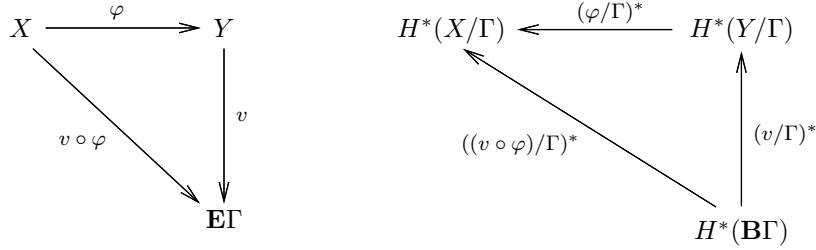

\begin{center}
  \begin{picture}(0,0)%
    \includegraphics{comm.pstex}%
  \end{picture}%
  \input{comm.pstex_t}%
 
\end{center}
\caption{Functoriality of characteristic classes.}
\label{fig:comm}
\end{figure}

Let $Y$ be another regular CW complex with a~free action of $\Gamma$,
and assume that $\varphi:X\ra Y$ is a~$\Gamma$-equivariant map. By
what is said above, there exists a~$\Gamma$-map $v:Y\ra{\bf E}\Gamma$. Hence,
in addition to the map $w:X\ra{\bf E}\Gamma$, we also have
a~composition map $v\circ\varphi$. Passing on to the quotient map and
then to the induced map on cohomology, we get yet another map
$((v\circ\varphi)/\Gamma)^*:H^*({\bf B}\Gamma)\ra
H^*(X/\Gamma)$. However, as we mentioned above, the map on the
cohomology algebras does not depend on the choice of the original map
to ${\bf E}\Gamma$. Thus, since $((v\circ\varphi)/\Gamma)^*=(\varphi/\Gamma)^*\circ(v/\Gamma)^*$, we have commuting diagrams, see Figure~\ref{fig:comm}, and therefore
	\[w(z,X)=(\varphi/\Gamma)^*(w(z,Y)),\]
where $z$ is an arbitrary element
of $H^*({\bf B}\Gamma)$.

In other words the characteristic classes associated to a~finite group
action are {\em natural}, or, as one sometimes says, {\em functorial}. 

We refer the reader to the wonderful book of tom~Dieck, \cite{tD}, for further details on equivariant maps and associated bundles. We also recommend the classical book of Milnor\&Stasheff, \cite{MS}, as an~excellent source for the theory of characteristic classes of vector bundles. The generalities on bundles, including principal bundles, can be found in~\cite{Stee}.

\subsection{$\zz$-spaces and the definition of Stiefel-Whitney classes.} 
\label{ssect:z2spaces} 


\nin Let now $X$ be a~$\zz$-space, i.e., a~CW complex equipped with 
a~fixed point free involution. Specifying $\Gamma=\zz$ in the
considerations above, we get a~map $w:X\ra S_a^{\infty}={\bf E}
\zz$. Furthermore, we have the induced quotient map
\[w/\zz:X/\zz\ra S_a^{\infty}/\zz=\rp^\infty={\bf B}\zz.\]
In this particular case, the induced $\zz$-algebra homomorphism
\[(w/\zz)^*:H^*(\rp^\infty;\zz)\ra H^*(X/\zz;\zz)
\]
is determined by very little data. Namely, let $z$ denote the nontrivial cohomology class in $H^1(\rp^\infty;\zz)$. Then $H^*(\rp^\infty;\zz) \simeq\zz[z]$ as a~graded $\zz$-algebra, with $z$ having degree~1. We denote the image $(w/\zz)^*(z)\in H^1(X/\zz;\zz)$ by $\sw(X)$. Obviously, the whole map $(w/\zz)^*$ is determined by the element~$\sw(X)$. 

This is the {\em Stiefel-Whitney class} of the $\zz$-space~X. Clearly,
$\sw^k(X)=(w/\zz)^*(z^k)$. Furthermore, by the general observation, if
$Y$ is another $\zz$-space, and $\varphi:X\ra Y$ is a~$\zz$-map, then
$(\varphi/\zz)^*(\sw(Y))=\sw(X)$.

As an example we can quickly compute $\sw(S_a^n)$, for an arbitrary nonnegative integer~$n$. First, for dimensional reasons, $\sw(S_a^0)=0$. So we assume $n\geq 1$. Next, we have $S_a^n/\zz= \rp^n$. The cohomology algebra $H^*(\rp^n;\zz)$ is generated by one element $\beta\in H^1(\rp^n;\zz)$, with a~single relation $\beta^{n+1}=0$. Finally, the standard inclusion map $\iota:S^n_a\hra S^\infty_a$ is $\zz$-equivariant, and induces another standard inclusion $\iota/\zz:\rp^n \hra\rp^\infty$. Identifying $\rp^n$  with the image of $\iota/\zz$, we can think of it as the $n$-skeleton of $\rp^\infty$. Thus the induced $\zz$-algebra homomorphism $(\iota/\zz)^*:H^*(\rp^\infty;\zz)\ra H^*(\rp^n;\zz)$ maps the canonical generator of $H^1(\rp^\infty;\zz)$ to $\beta$, and hence we can conclude that $\sw(S_a^n)=\beta$.

\section{Properties of Stiefel-Whitney classes.}

\subsection{Borsuk-Ulam theorem, index, and coindex.} \label{ssect:index} 


\nin The Stiefel-Whitney classes can be used to determine the nonexistence of certain $\zz$-maps. The following theorem is an example of such situation.

\begin{thm} {\em (Borsuk-Ulam).}\label{thm:BU}

\nin Let $n$ and $m$ be nonnegative integers. If there exists a~$\zz$-map $\varphi:S^n_a\ra S^m_a$, then $n\leq m$. 
\end{thm}

\pr Choose representations for the cohomology algebras $H^*(\rp^n;\zz)=\zz[\alpha]$, and $H^*(\rp^m;\zz)=\zz[\beta]$, with the only relations on the generators being $\alpha^{n+1}=0$, and $\beta^{m+1}=0$. Since the Stiefel-Whitney classes are functorial, we get $(\varphi/\zz)^*(\sw(S^m_a))=\sw(S^n_a)$.

On the other hand, by the computation in the subsection~\ref{ssect:z2spaces}, we have $\sw(S^n_a)=\alpha$, and $\sw(S^m_a)=\beta$. So $(\varphi/\zz)^*(\beta)=\alpha$, and hence $\alpha^{m+1}=(\varphi/\zz)^*(\beta)^{m+1}=(\varphi/\zz)^*(\beta^{m+1})=0$. Since $\alpha^{n+1}=0$ is the only relation on $\alpha$, this yields the desired inequality~$m\geq n$.
\qed

\vskip5pt

The Borsuk-Ulam Theorem makes the following terminology useful for formulating further obstructions to maps between $\zz$-spaces.


\begin{df} \label{df:index}
Let $X$ be a $\zz$-space. 

\nin $\bu$ The {\bf index} of $X$, denoted $\ind X$, is the minimal integer $n$, for which there exists a~$\zz$-map from $X$ to $S_a^n$. 

\nin $\bu$ The {\bf coindex} of $X$, denoted $\coind X$, is the maximal integer $n$, for which there exists a~$\zz$-map from $S_a^n$ to $X$.   
\end{df}

Assume that we have two $\zz$-spaces $X$ and $Y$, and that
$\gamma:X\ra Y$ is a~$\zz$-map, then, we have the inequality $\coind
X\leq\ind Y$. Indeed, if there exists $\zz$-maps $\varphi:S_a^n\ra X$,
and $\psi:Y\ra S_a^m$, then the composition
\[S_a^n\stackrel{\varphi}{\lra}X\stackrel{\gamma}{\lra}Y
\stackrel{\psi}{\lra}S_a^m\]
yields a~$\zz$-map between two spheres with antipodal actions, hence, by the Borsuk-Ulam Theorem, we can conclude that $n\leq m$. In particular, taking $Y=X$, and $\varphi=\id$, we get the inequality $\coind X\leq\ind X$, for an arbitrary $\zz$-space.

\subsection{Higher connectivity and Stiefel-Whitney classes.} \label{ssect:kconn} 


\nin Many results giving topological obstructions to graph colorings had the $k$-connectivity of some space as the crucial assumption. We notice here an important connection between this condition and non-nullity of powers of Stiefel-Whitney classes.

First, it is trivial, that if $X$ is a~non-empty $\zz$-space, then one can equi\-variantly map $S_a^0$ to $X$. It is possible to extend this construction inductively to an~arbitrary $\zz$-space.

\begin{prop} \label{pr:conn}
Let $X$ and $Y$ be two simplicial complexes with a~free $\zz$-action, such that for some $k\geq 0$, we have $\dim X\leq k$, and  $Y$ is $(k-1)$-connected. Assume further that we have a~$\zz$-map $\psi:X^{(d)}\ra Y$, for some $d\geq -1$. Then, there exists a~$\zz$-map $\varphi:X\ra Y$, such that $\varphi$ extends~$\psi$.
\end{prop}

Please note the following convention used in the formulation of Proposition~\ref{pr:conn}: $d=-1$ means we have no map $\psi$ (in other words, $X^{-1}=\emptyset$), hence no additional conditions on the map~$\varphi$.

\vskip5pt

\pr Choose a~$\zz$-invariant simplicial structure on $X$. We construct $\varphi$ inductively on $i$-skeleton of $X$, for $i\geq d+1$. If $d=-1$, we start by defining $\varphi$ on the 0-skeleton as follows: for each orbit $\{a,b\}$ consisting of two vertices of $X$, simply map $a$ to an arbitrary point $y\in Y$, and then map $b$ to $\gamma(y)$, where $\gamma$ is the free involution of~$X$. 

Assume now that $\varphi$ is defined on the $(i-1)$-skeleton of $X$, and extend the construction to the $i$-skeleton as follows. Let $(\sigma,\tau)$ be a~pair of $i$-dimensional simplices of $X$, such that $\gamma\sigma=\tau$. The boundary $\partial\sigma$ is a~$(i-1)$-dimensional sphere. By our assumptions $i-1\leq \dim X-1\leq k-1$, hence the restriction of $\varphi$ to $\partial\sigma$ extends to $\sigma$.
Finally, we extend $\varphi$ to the second simplex $\tau$ by applying the involution $\gamma$: $\varphi|_{\tau}:=(\varphi|_{\sigma})\circ\gamma$.
 \qed

\begin{crl} \label{crl:conn}
Let $X$ be a~$\zz$-space, and assume $X$ is $(k-1)$-connected, for some $k\geq 0$. Then there exists a~$\zz$-map $\varphi:S_a^k\ra X$. In particular, we have $\sw^k(X)\neq 0$.
\end{crl}

\pr $S_a^k$ is $k$-dimensional, hence the statement follows immediately from Proposition~\ref{pr:conn}. To see that $\sw^k(X)\neq 0$, recall that, since the Stiefel-Whitney classes are functorial, we have $(\varphi/\zz)^* (\sw^k(X)) =\sw^k(S_a^k)$, and the latter has been verified to be nontrivial.
\qed

\vskip5pt

The Corollary~\ref{crl:conn} explains the rule of thumb that, whenever dealing with $\zz$-spaces, the condition of $k$-connectivity can be replaced by the weaker condition that the $(k+1)$-th power of the appropriate Stiefel-Whitney class is different from~0.

\subsection{Combinatorial construction of Stiefel-Whitney classes.} \label{ssect:swcomb} 


\nin Let us describe how the construction used in the proof of Proposition~\ref{pr:conn} can be employed to obtain an~explicit combinatorial description of the Stiefel-Whitney classes.

Let $X$ be a~regular CW complex and a~$\zz$-space, and denote the fixed point free involution on $X$ by~$\gamma$. As mentioned above, one can choose a~simplicial structure on $X$, such that $\gamma$ is a~simplicial map. We define a~$\zz$-map $\varphi:X\ra S_a^\infty$ following the recipe above.

Take the standard $\zz$-equivariant cell decomposition of $S_a^\infty$ with two antipodal cells in each dimension. Divide $X^{(0)}$, the set of the vertices of $X$, into two disjoint sets $X^{(0)}=A\cup B$, such that every orbit of the $\zz$-action contains exactly one element from $A$ and one element from~$B$.  Let $\{a,b\}$ be the 0-skeleton of $S_a^\infty$, and map all the points in $A$ to $a$, and all the points in $B$ to $b$. Call the edges having one vertex in $A$, and one vertex in $B$, {\em multicolored}, and the edges connecting two vertices in $A$, resp.\ two vertices in $B$, $A$-internal, resp.\ $B$-internal.

Let $\{e_1,e_2\}$ be the 1-skeleton of $S_a^\infty$. One can then
extend $\varphi$ to the 1-skeleton as follows. Map the $A$-internal
edges to $a$, map the $B$-internal edges to $b$. Note that the
multicolored edges form $\zz$-orbits, 2 edges in every orbit. For each
such orbit, map one of the edges to $e_1$ (there is some arbitrary
choice involved here), and map the other one to~$e_2$.

Since the $\zz$-action on the space $X$ is free, the generators of the cochain complex $C^*(X/\zz;\zz)$ can be indexed with the orbits of simplices. For an~arbitrary simplex $\delta$ we denote by $\tau_\delta$ the generator corresponding to the orbit of $\delta$; in particular, $\tau_{\gamma(\delta)}=\tau_\delta$. 

The induced quotient cell decomposition of $\rp^\infty$ is the standard one, with one cell in each dimension. The cochain $z^*$, corresponding to the unique edge of $\rp^\infty$, is the generator (and the only nontrivial element) of $H^1(\rp^\infty;\zz)$. Its image under $(\varphi/\zz)^*$ is simply the sum of all orbits of the multicolored edges:
\begin{equation} \label{eq:multedge}
 \sw(X)=(\varphi/\zz)^*(z^*)=\sum_{\text{multicolored }e}\tau_e,
\end{equation}
where the sum is taken over representatives of $\zz$-orbits of multicolored edges, one representative per orbit.

To describe the powers of the Stiefel-Whitney classes, $\sw^k(X)$, we need to recall how the cohomology multiplication is done simplicially. In fact, to evaluate $\sw^k(X)$ on a~$k$-simplex $(v_0,v_1,\dots,v_k)$, we need to evaluate $\sw(X)$ on each of the edges $(v_i,v_{i+1})$, for $i=0,\dots,k-1$, and then multiply the results. Thus, the only $k$-simplices, on which the power $\sw^k(X)$ evaluates nontrivially, are those whose ordered set of vertices has alternating elements from $A$ and from $B$. We call these simplices {\em multicolored.} We summarize
\begin{equation} \label{eq:multsimp}
 \sw^k(X)=\sum_{\text{multicolored }\sigma}\tau_\sigma,
\end{equation}
where the sum is taken over representatives of $\zz$-orbits of multicolored $k$-dimensional simplices, one representative per orbit.

\section{First applications of Stiefel-Whitney classes to lower bounds of chromatic numbers of graphs.}

\subsection{Complexes of complete multipartite subgraphs.} \label{ssect:homkn} 


\nin We have seen in the subsection~\ref{ss3.1} that the complex $\thom(K_2,G)$ is simply homotopy equivalent to the neighbourhood complex of $G$. In particular, $\thom(K_2,K_n)$ is simply homotopy equivalent to~$S^{n-2}$. The following proposition provides us with a~more complete information.

\begin{prop}\label{pr:k2kn} $\,$ {\rm (\cite[Proposition 4.3]{BK03b}).}

\nin {\rm (a)} $\thom(K_2,K_{n+1})$ is isomorphic as a~cell complex to the~boundary complex of the Minkowski sum $\Delta^n+(-\Delta^n)$. 

\nin {\rm (b)}  The $\zz$-action on $\thom(K_2,K_{n+1})$, induced by the flip action of $\zz$ on $K_2$, corresponds under this isomorphism to the central symmetry.
\end{prop} 

\begin{figure}[hbt]
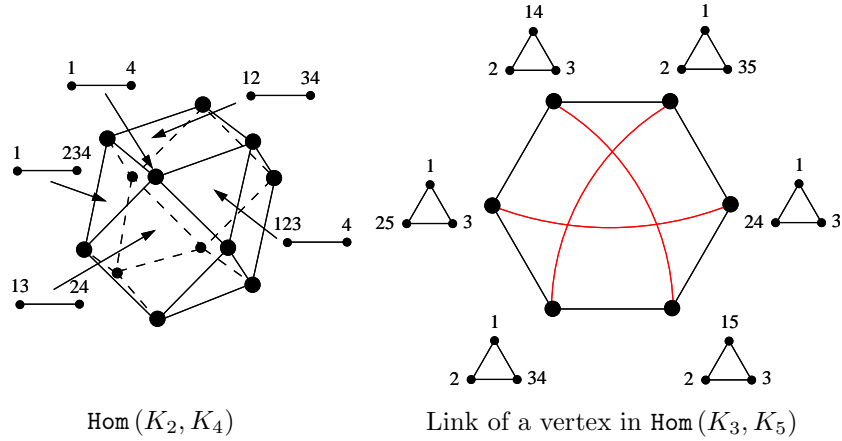

\begin{center}
  \begin{picture}(0,0)%
    \includegraphics{hom3a.pstex}%
  \end{picture}%
  \input{hom3a.pstex_t}%
 
\end{center}
\caption{Complexes of graph homomorphisms between complete graphs.}
\label{fig:hom3a}
\end{figure}

Proposition~\ref{pr:k2kn} is illustrated with Figure~\ref{fig:hom3a}. Perhaps the easiest way to Proposition~\ref{pr:k2kn} is by means of the following notion.

\begin{df} \label{df:delprod}
Let $X_1,\dots,X_t$ be a family of simplicial complexes with isomorphic sets of vertices. The {\bf deleted product} of this family is the subcomplex of the direct product of $X_1,\dots,X_t$, consisting of all cells $\tau_1\times\dots\times\tau_t$, satisfying $\tau_i\cap\tau_j=\emptyset$, for any $i\neq j$.
\end{df}
 
Clearly, $\thom(K_m,K_n)$ can be viewed as a~deleted product of $m$ copies of $(n-1)$-dimensional simplices, see e.g.,~\cite{Ma1}. In this context Proposition~\ref{pr:k2kn} is well-known, probably due to van~Kampen.

Also, for $m\geq 3$, $\thom(K_m,G)$ can be thought of as complexes consisting of all complete $m$-partite subgraphs of $G$. Even in the case $G=K_n$, it seems complicated to understand $\thom(K_m,G)$ up to homeomorphism. However, we still obtain a~good description of the homotopy type.
  
\begin{thm} \label{thm:kmkn} {\rm (\cite[Proposition 4.5]{BK03b}).}

\nin The cell complex $\thom(K_m,K_n)$ is homotopy equivalent to a~wedge of $(n-m)$-dimensional spheres.
\end{thm}

Introducing a~new piece of notation, let us say that the complex $\thom(K_m,K_n)$ is homotopy equivalent to a~wedge of $f(m,n)$ spheres.
Let $S(-,-)$ denote the {\em Stirling numbers} of the second kind, and $SF_k(x)=\sum_{n\geq k}S(n,k)x^n$ denote the generating function (in the first variable) for these numbers. It is well-known, see for example \cite[p.\ 34]{St1}, that 
$$SF_k(x)=x^k/((1-x)(1-2x)\dots(1-kx)).$$ 
For $m\geq 1$, let $F_m(x)=\sum_{n\geq 1}f(m,n)x^n$ be the generating
function (in the second variable) for the number of the spheres. Clearly, $F_1(x)=0$, and $F_2(x)=x^2/(1-x)$.

\begin{prop} {\rm (\cite[Proposition 4.6]{BK03b}).}

\nin 
The numbers $f(m,n)$ satisfy the following recurrence relation
\begin{equation}\label{eqrecf}
f(m,n)=mf(m-1,n-1)+(m-1)f(m,n-1),
\end{equation}
for $n>m\geq 2$; with the boundary values $f(n,n)=n!-1$, $f(1,n)=0$
for $n\geq 1$, and $f(m,n)=0$ for $m>n$.

Then, the generating function $F_m(x)$ is given by the equation:
\begin{equation}\label{eqngenff}
F_m(x)=(m!\cdot x \cdot SF_{m-1}(x)-x^m)/(1+x).
\end{equation}

As a~consequence, the following non-recursive formulae are valid:
\begin{equation}\label{eqnrecf2}
f(m,n)=(-1)^{m+n+1}+m!(-1)^{n}\sum_{k=m}^{n}(-1)^k S(k-1,m-1),
\end{equation}
and 
\begin{equation}\label{eqnrecf}
f(m,n)=\sum_{k=1}^{m-1}(-1)^{m+k+1}{\binom{m}{k+1}} k^n,
\end{equation}
for $n\geq m\geq 1$.
\end{prop}

In particular, for small values of $m$, we obtain the following explicit formulae: $f(2,n)=1$, for $n\geq 2$, $f(3,n)=2^n-3$, for $n\geq 3$, $f(4,n)=3^n-4\cdot 2^n+6$, for $n\geq 4$, $f(5,n)=4^n-5\cdot 3^n+10\cdot 2^n-10$, for $n\geq 5$.

\subsection{Stiefel-Whitney classes and test graphs.} \label{ssect:swtest} 


\nin One connection between the non-nullity of the powers of Stiefel-Whitney characteristic classes and the lower bounds for graph colorings is provided by the following general observation.

\begin{thm} \label{thm:swlb}
Let $G$ be a graph without loops, and let $T$ be a~graph with
$\zz$-action which flips some edge in $T$. If, for some integers
$k\geq 0$, $m\geq 1$, $\sw^k(\thom(T,G))\neq 0$, but
$\sw^{k}(\thom(T,K_m))=0$, then $\chi(G)\geq m+1$.
\end{thm}

Since this statement is crucial for all our applications, we provide here a~short argument.

\vskip5pt

\nin {\bf Proof of Theorem~\ref{thm:swlb}.} We know that, under the assumptions of the theorem, $\thom(T,H)$ is a~$\zz$-space for any loopfree graph $H$. Assume now that the graph $G$ is $m$-colorable, i.e., there exists a~homomorphism $\varphi:G\ra K_m$. It induces a~$\zz$-map $\varphi^{T}:\thom(T,G)\ra\thom(T,K_m)$. Since the
Stiefel-Whitney classes are functorial and $\sw^k(\thom(T,K_m))=0$, the existence of the map $\varphi^{T}$ implies that $\sw^k(\thom(T,G))=0$, which is a~contradiction to the assumption of the theorem. \qed

\vskip5pt

We can now use Theorem~\ref{thm:kmkn} to give lower bounds for
chromatic numbers of graphs in terms of Stiefel-Whitney classes of
complexes of graph homomorphisms from complete graphs.

\begin{thm}\label{thm:cothm}
Let $G$ be a graph, and let $n,k\in\dz$, such that $n\geq 2$, $k\geq
-1$. If $\sw^k(\thom(K_n,G))\neq 0$, then $\chi(G)\geq k+n$.
\end{thm}

\pr Indeed, substituting $T=K_n$, and $m=k+n-1$, in the
Theorem~\ref{thm:swlb}, all we need to do is to see that
$\sw^k(\thom(K_n,K_{k+n-1}))=0$. By Theorem~\ref{thm:kmkn},
$\thom(K_n,K_{k+n-1})$ is homotopy equivalent to a~wedge of
$(k-1)$-dimensional spheres. Hence, by dimensional reasons we conclude that 
$\sw^k(\thom(K_n,K_{k+n-1}))=0$. \qed

\lecture{The spectral sequence approach.} 

\section{$\thomp$-construction.}

\subsection{Various definitions.} \label{ss6.1} 


\nin We shall now define a complex $\thomp(T,G)$ which is related to $\thom(T,G)$. It is easier to compute various algebro-topological invariants for this complex, however, it also has less bearing on our main problem: computation of the lower bounds for chromatic numbers. We shall then connect the $\thom$- and $\thomp$-constructions by means of a~spectral sequence.

\subsubsection{A subcomplex of a total join.} \label{sss6.1.1} $\,$

\vskip5pt

\nin Let $T$ and $G$ be arbitrary graphs. We shall define $\thomp(T,G)$ analogously to $\thom(T,G)$, replacing the direct product with the join. We note here that whenever talking about $\thomp(T,G)$ we always assume that the graph $T$ is finite.

Let, as before, $\Delta^{V(G)}$ be a~simplex whose set of vertices is
$V(G)$. Let $J(T,G)$ denote the join $\join{x\in V(T)} \Delta^{V(G)}$,
i.e., the copies of $\Delta^{V(G)}$ are indexed by vertices of~$T$.
A~cell (simplex) in $J(T,G)$ is a~join of (possibly empty) simplices
$\join{x\in V(T)}\sigma_x$, the dimension of this simplex is
$\sum_{x\in V(T)} (\dim\sigma_x+1)-1$. Remark that this number is
finite, since we assumed that $T$ is finite. Here we use the usual
convention that $\dim\emptyset=-1$.

\begin{df} \label{df:homp}
For arbitrary graphs $T$ and $G$,
$\thomp(T,G)$ is the simplicial subcomplex of $J(T,G)$ defined by the following condition: $\sigma=\join{x\in V(T)}\sigma_x\in\thomp(T,G)$ if and only if for any $x,y\in V(T)$, if $(x,y)\in E(T)$, and both $\sigma_x$ and $\sigma_y$ are nonempty, then $(\sigma_x,\sigma_y)$ is a~complete bipartite subgraph of~$G$.
\end{df}

The intuition behind this definition is that we relax the conditions of the $\thom$ case by allowing some of the "coloring lists" to be empty.
One can think of $\thomp(T,G)$ as a~simplicial structure imposed on the set of all {\it partial graph homomorphisms} from $T$ to~$G$, i.e., graph homomorphisms from an induced subgraph of $T$ to the graph~$G$.

In analogy with the $\thom$ case, we can describe the simplices of $\thomp(T,G)$ directly: they are indexed by all $\eta:V(T)\ra 2^{V(G)}$ satisfying the same condition as in Definition~\ref{df:hom}. The closure of $\eta$ is also defined identical to how it was defined for $\thom$. So the only difference is that $\eta(x)$ is allowed to be an empty set, for $x\in V(T)$.

\subsubsection{A link of a~vertex in an auxiliary $\thom$-complex.} \label{sss6.1.2} $\,$

\vskip5pt

\nin The following construction is the graph analog of the topological coning.

\begin{df} \label{df:plus}
For an~arbitrary graph $G$, let $G_+$ be the graph obtained from $G$ by
adding an extra vertex $a$, called the apex vertex, and connecting it
by edges to all the vertices of $G_+$ including $a$ itself, i.e.,
$V(G_+)=V(G)\cup\{a\}$, and $E(G_+)=E(G)\cup\{(x,a),(a,x)\,|\,x\in
V(G_+)\}$.
\end{df}

We note that, for an arbitrary polyhedral complex $K$, such that all
faces of $K$ are direct products of simplices, and a~vertex $x$
of~$K$, the {\it link} of $x$, $\link_K(x)$, is a~simplicial complex.
It follows from the fact that a~link of any vertex in a~hypercube is
a~simplex, and the identity $\link_{(A\times
B)}(v,w)=\link_A(v)*\link_B(w)$, for arbitrary polyhedral complexes $A$
and~$B$.

We are now ready to formulate another definition, which is equivalent to Definition~\ref{df:homp}.

\begin{df} \label{df:homp2}
For arbitrary graphs $T$ and $G$, the simplicial complex $\thomp(T,G)$ is defined to be the link in $\thom(T,G_+)$ of the specific graph homomorphism $\alpha$, which maps all vertices of $T$ to the apex vertex of $G_+$. 
\end{df}

In short: $\thom_+(T,G)=\link_{\thom(T,G_+)}(\alpha)$.

The equivalence of the definitions follows essentially from the
following bijection: let $\eta\in\cf(\thom(T,G_+))_{>\alpha}$, and set
$\ti\eta(v):=\eta(v)\sm\{a\}$, for any $v\in V(T)$. Clearly,
$\ti\eta\in\cf(\thomp(T,G))$, and it is easily checked that this
bijection produces an isomorphism of simplicial complexes.

\subsubsection{Functorial properties of the $\thomp$-construction.} \label{sss6.1.3} $\,$

\vskip5pt

\nin Just like in the case of the $\thom(-,-)$-construction, $\thomp(T,-)$ is a~covariant functor from {\bf~Graphs} to {\bf~Top}. For two arbitrary graphs $G$ and $K$, and a~graph homomorphism $\varphi$ from $G$ to $K$, we have an~induced simplicial map $\varphi^T:\thomp(T,G)\ra\thomp(T,K)$.

Again, as in the case of $\thom(-,-)$, the situation is somewhat more complicated with the functoriality in the first argument. Let $T,G$, and $K$, be three arbitrary graphs. This time, for a~graph homomorphism $\psi$ from $T$ to $G$ to induce a~topological map from $\thomp(G,K)$ to $\thomp(T,K)$, we must require that $\psi$ is surjective on the vertices. We can define the topological map $\psi_K$ in the same way as for the $\thom(-,-)$ case, but if $\psi$ is not surjective on the vertices, then we may end up mapping a~non-empty cell to an~empty one. If, in addition, we want a~simplicial map $\psi_K:\thomp(G,K)\ra\thomp(T,K)$, then, as before in the subsection~\ref{ss:231}, we must require that $\psi$ is injective, hence bijective on the vertices. 

In particular, we still have that the group $\aut(T)\times\aut(G)$ acts on the complex $\thomp(T,G)$ simplicially. The difference is that we do not have the freeness as easily as we had in the $\thom(-,-)$ case. For example, for an~involution $\gamma$ of $T$ to induce a~free action $\gamma_G$ on $\thomp(T,G)$ we need to require that all orbits of $\gamma$ on $V(T)$ are of cardinality~2, and that the vertices in the same orbit are connected by an~edge. For instance, the action of $\zz$ on $\thomp(K_2,G)$ is free, whereas the reflection $\zz$-action on $\thomp(C_{2r+1},G)$ is not.

\subsection{Connection to independence complexes.} \label{ss6.2} 


\nin The following is a standard construction in topological combinatorics, see e.g., \cite{K2,Me1}.

\begin{df} \label{df:indcomp}
For an arbitrary graph $G$, the {\bf independence complex} of $G$, $\ind(G)$, is the simplicial complex, whose set of vertices is $V(G)$, and simplices are all the independent sets (anticliques) of~$G$.
\end{df}

Before we can make use of the $\ind(-)$-construction in our context, we need more graph terminology.

\begin{df} \label{df:strcomp}
For an arbitrary graph $G$, the {\bf strong complement} $\cmp G$ is
defined by $V(\cmp G)=V(G)$, and $E(\cmp G)=V(G)\times V(G)\sm E(G)$.
\end{df} 

For example, $\cmp K_n$ is the disjoint union of $n$ loops.

\begin{df} \label{df:dirprod}
For arbitrary graphs $G$ and $H$, the {\bf direct product} $G\times H$ is defined by: $V(G\times H)=V(G)\times V(H)$, and $E(G\times H)= \{((x,y),(x',y'))\,|\,(x,x')\in E(G),(y,y')\in E(H)\}$.
\end{df}

For example, $K_2\times K_2$ is a disjoint union of two copies of $K_2$, whereas $G\times\cmp K_1$ is isomorphic to $G$ for an arbitrary graph~$G$.

\begin{figure}[hbt]
\begin{center}
  \begin{picture}(0,0)%
    \includegraphics{ex+.pstex}%
  \end{picture}%
  \input{ex+.pstex_t}%
 
\end{center}
\caption{The $+$-construction.}
\label{fig:ex+}
\end{figure}

Sometimes, it can be convenient to view $\thomp(G,H)$ as the independence
complex of a~certain graph.

\begin{prop} \label{pr:homp}
{\rm (\cite[Proposition 3.2]{BK03c}).}

\nin For arbitrary graphs $T$ and $G$,
$\thomp(T,G)$ is isomorphic to $\ind(T\times\cmp G)$. 
\end{prop}

Specializing Proposition~\ref{pr:homp} to $G=K_n$, and taking into account $\cmp K_n= \coprod_{i=1}^n \cmp K_1$ (observed above), and the fact that for arbitrary graphs $G_1$ and $G_2$ we have
        \[\ind(G_1\coprod G_2)=\ind(G_1)*\ind(G_2),
\]
 we obtain the following corollary.

\begin{crl} \label{crl:hompkn}
For an~arbitrary graph $T$,
$\thomp(T,K_n)$ is isomorphic to the $n$-fold join $\ind(T)^{*n}$.
\end{crl}

When $G$ is loopfree, the~dimension of the simplicial complex $\thomp(T,G)$ (unlike that of $\thom(T,G)$) is easy to find, once the size of the maximal independent set of $G$ is computed.
 
\begin{prop} \label{crl:dimhomp}
For an arbitrary graph $T$, and an arbitrary loopfree graph~$G$, we have
\[
\dim(\thomp(T,G))=|V(G)|\cdot(\dim(\ind(T))+1)-1.
\]
\end{prop}

\pr Indeed, let $s=\dim(\ind(T))+1$ be the size of the maximal independent set in $T$. Since $G$ is loopfree, every vertex of $G$ occurs in at most $s$ of the sets $\eta(x)$, for $x\in V(T)$. On the other hand, we can choose an~independent set $S\subseteq V(T)$, such that $|S|=s$, and then
assign 
        \[\eta(x)=\begin{cases} 
        V(G),&\text{ for }x\in S;\\
        \emptyset,& \text{ otherwise.}
        \end{cases}
\]
This gives a simplex of dimension $|V(G)|\cdot(\dim(\ind(T))+1)-1$.\qed

\vskip5pt

For example, $\dim(\thomp(C_{2r+1},K_n))=n\cdot((r-1)+1)-1=nr-1$.

\subsection{The support map.} \label{ss6.3} 


\nin For any topological space $X$ and a set $I$, there is the standard support map from the join of $I$ copies of $X$ to the appropriate simplex
        \[\supp: *_I X\longrightarrow \Delta^I,
\]
which "forgets" the coordinates in $X$.

Specializing to our situation, for arbitrary graphs $T$ and $G$, we get the restriction map $\supp:\thomp(T,G)\ra\Delta^{V(T)}$. Explicitly, for each simplex of $\thomp(T,G)$, $\eta:V(T)\ra 2^{V(G)}$, the support of $\eta$ is given by $\supp\eta=V(T)\setminus\eta^{-1}(\emptyset)$. 

An~important property of the support map is that the preimage of the
barycenter of $\Delta^{V(T)}$ is homeomorphic to $\thom(T,G)$. This is
the crucial step in setting up a~useful spectral sequence. The
assumption that $T$ is finite is crucial at this point, since an
infinite simplex does not have a~barycenter.

\begin{figure}[hbt]
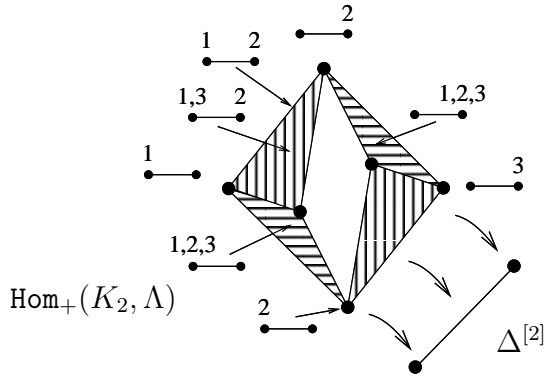

\begin{center}
  \begin{picture}(0,0)%
    \includegraphics{supp.pstex}%
  \end{picture}%
  \input{supp.pstex_t}%
 
\end{center}
\caption{The support map from $\thomp(K_2,\Lambda)$ to $\Delta^{[2]}$.}
\label{fig:supp}
\end{figure}

An~alternative concise way to phrase the definition of $\supp$ is to consider the map $t^T:\thomp(T,G)\ra\thomp(T,\cmp K_1)\simeq\Delta^{V(T)}$ induced by the homomorphism $t:G\ra\cmp K_1$. Then, for each $\eta\in\thomp(T,G)$ we have $\supp\eta=t^T(\eta)$, where the simplices in $\Delta^{V(T)}$ are identified with the finite subsets of~$V(T)$.

\section{Spectral sequence generalities.} \label{sect:spec} 

\nin Spectral sequences constitute an important tool of topological combinatorics in general. They have also been proved invaluable in the solution of the Lov\'asz Conjecture. Taking into account the format of this article, we only give a short introduction here, aimed at setting up the notations and, hopefully, helping the intuition. We refer the interested reader to the excellent existing sources, see e.g., \cite{FF,McC}.

\subsection{Cochain complexes and their cohomology.} \label{ssect:chcomp} 


\nin Recall that a {\em cochain complex} is a sequence
\[
\cc=\dots \stackrel{d^{i-2}}{\lra}C^{i-1}\stackrel{d^{i-1}}{\lra} C^i \stackrel{d^i}{\lra}C^{i+1} \stackrel{d^{i+1}}{\lra}\dots,
\]
where $C^i$'s are $\car$-modules, $d^i$'s are $\car$-module homomorphisms (called differentials) satisfying $d^{i+1}\circ d^i=0$, and $\car$ is a~commutative ring with a unit. We shall use the notation $\cc=(C^*,d^*)$. For all our purposes it is enough in the continuation to restrict one's attention to the cases $\car=\dz$, and $\car=\zz$.

We can associate a cochain complex $C^*(X)$ to a cell complex $X$ in the standard way: $C^i(X)$ is taken to be a~free $\car$-module with the generators indexed by the $i$-dimensional cells (the module of $\car$-valued functionals on cells of $X$), and the differential maps are given by the corresponding coboundary maps. Sometimes this particular cochain complex is called a~{\em cellular cochain complex}.

Given two cochain complexes $\cc_1=(C^*_1,d^*_1)$ and $\cc_2=(C^*_2,d^*_2)$, 
a~{\em cochain complex homomorphism} $\varphi:\cc_1\ra\cc_2$ (also called a~{\em cochain complex map}) is a~collection of $\car$-module homomorphisms $\varphi_i:C_1^i\ra C_2^i$, for all integers $i$, such that the following diagram commutes
\begin{equation}
\begin{CD}
\label{eq:chainmap}
C_1^i @>d_1^i>> C_1^{i+1} \\ 
@V\varphi_i VV  @VV\varphi_{i+1} V\\
C_2^i @>d_2^i>>C_2^{i+1} \\
\end{CD}
\end{equation}

For each choice of $\car$, cochain complexes together with cochain complex homomorphisms form a~category.
 
Associated with a cochain complex, one has the {\em cohomology groups }
\[  H^i(\cc)=\ker d^i/\im d^{i-1}.
\]
In our cases $H^i(\cc)$ is either an abelian group or a vector space over~$\zz$. 

Given two cochain complexes $\cc_1=(C^*_1,d^*_1)$ and
$\cc_2=(C^*_2,d^*_2)$, and a~cochain complex homomorphism
$\varphi:\cc_1\ra\cc_2$, since $\varphi_i$'s commute with the
corresponding differentials, $\varphi$ induces a~map on the cohomology
groups $\varphi^*_i:H^i(\cc_1)\ra H^i(\cc_2)$.

The above facts mix well with the cellular structure. First, for a~cell complex $X$, the cellular cohomology groups of $X$ are by definition isomorphic to the cohomology groups of the associated cochain complex $C^*(X)$. Second, for two cell complexes $X$ and $Y$, a~cellular map $\varphi:X\ra Y$ induces a~cochain complex map between associated cochain complexes (but in the opposite direction!), and hence a~map between corresponding cohomology groups.

\subsection{Filtrations.} \label{ssect:filtr} 


\nin In concrete situations it can be difficult to compute the 
cohomology groups $H^i(\cc)$ without auxiliary constructions. The idea
behind spectral sequences is to break up this large task into smaller
subtasks, with the formal machinery to help the bookkeeping. This
"break up" is usually phrased in terms of a~filtration.

A {\em cochain subcomplex} of $\cc$ is a sequence
\[
\wti\cc=\dots \stackrel{d^{i-2}}{\lra}\wti C^{i-1}\stackrel{d^{i-1}}{\lra} \wti C^i \stackrel{d^i}{\lra}\wti C^{i+1} \stackrel{d^{i+1}}{\lra}\dots,
\]
where $\wti C^i$ is an $\car$-submodule of $C^i$, and the differentials are restrictions of those in $\cc$. We shall simply write $\cc\supseteq\wti\cc$.
In this situation, one can form the quotient cochain complex
\[
\cc/\wti\cc=\dots \stackrel{d^{i-2}}{\lra}C^{i-1}/\wti C^{i-1} \stackrel{d^{i-1}}{\lra} C^i/\wti C^i \stackrel{d^i}{\lra} 
C^{i+1}/\wti C^{i+1} \stackrel{d^{i+1}}{\lra}\dots.
\]
The cohomology groups of this complex, $H^*(\cc/\wti\cc)$, are usually denoted $H^*(\cc,\wti\cc)$, and are called the {\em relative cohomology groups.} 

If $X$ is a~cell complex, and $Y$ its cell subcomplex, then the cellular cochain complex of $Y$ is a~cochain subcomplex of the cellular cochain complex of~$X$. The corresponding cohomology groups of the quotient cochain complex are precisely the relative cohomology groups of the pair of topological spaces $(X,Y)$.

\begin{df} \label{df:filtr}
A (finite) {\bf filtration} on a cochain complex $\cc$ is a nested sequence of cochain complexes 
\[\cc_j=\dots \stackrel{d^{i-2}}{\lra}C_j^{i-1}\stackrel{d^{i-1}}{\lra} C^i_j \stackrel{d^i}{\lra}C^{i+1}_j \stackrel{d^{i+1}}{\lra}\dots,
\]
for $j=0,1,2,\dots,t$, such that $\cc=\cc_t\supseteq\cc_{t-1}\supseteq\dots\supseteq\cc_0$ (that is why we suppressed the lower index in the differential).
\end{df}

In general, infinite filtrations can be considered, but in this article we limit our considerations to the finite ones. Given a~filtration $\cc=\cc_t\supseteq\cc_{t-1}\supseteq\dots\supseteq\cc_0$, we set $\cc_{-1}=0$, for the convenience of notations.

There are many standard filtrations of cochain complexes. For example, if a~pure cochain complex is bounded, say $C^i=0$, for $i<0$, or $i>t$, then, the standard skeleton filtration is defined as follows:
\[ C_j^i=\begin{cases}
C^i,& \text{ if } i\leq j; \\
0,& \text{ otherwise.}
\end{cases}
\]
This filtration is not very interesting though, since computing the cohomology groups with its help is canonically equivalent with computing the cohomology groups from the cochain complex directly.

For a cell complex $X$, a~classical way to define a~filtration on its cellular cochain complex, is to choose a~cell filtration on $X$, i.e., a~sequence of cell subcomplexes $X=X_t\supseteq X_{t-1}\supseteq\dots \supseteq X_0$ (again for the convenience of notations, we set $X_{-1}=\emptyset$). As mentioned above, the corresponding cellular cochain complexes form a~sequence of nested subcomplexes.

If the cell complex $X$ is finite dimensional, then, taking $X_i$ to be the $i$-th skeleton of $X$, we recover the standard skeleton filtration on $C^*(X)$, which explains the name of this filtration. 

A much more interesting situation is the following. 
\begin{df} \label{df:pullback}
Assume that we have a~cell map $\varphi:X\ra Y$ and a~filtration $Y=Y_t\supseteq Y_{t-1} \supseteq\dots \supseteq Y_0$. Define a~filtration on $X$ as follows:
$X_i:=\varphi^{-1}(Y_i)$, for $i=0,\dots,t$. This filtration on $X$ is called the {\bf pullback} of the filtration on $Y$ {\bf along $\varphi$.}
\end{df}
In the case when the filtration on $Y$ is simply the skeleton filtration, the corresponding pullback filtration on $X$ is called the {\em Serre filtration}. We use the same name for the corresponding filtration on the cellular cochain complex of~$X$.

\subsection{Spectral sequence terminology.} \label{ssect:specgen} 


\nin Once we have fixed a filtration $\cc=\cc_t\supseteq\cc_{t-1}\supseteq \dots\supseteq\cc_0$, $\cc_i=(C_i^*,d^*)$, on a~cochain complex, we can proceed to compute its cohomology groups by studying auxiliary algebraic gadgets derived from the filtration. 

Rather than studying the 1-dimensional cochain complex directly, we study a~sequence of 2-dimensional tableaux $E_n^{*,*}$, $n=0,1,2,\dots$. Our cochain complex had the usual differential, going one up in degree, which one can express symbolically by writing $d^1:C^*\ra C^{*+1}$. Instead, each tableau $E_n^{*,*}$ is equipped with a differential going almost diagonally, $d_n^{*,*}:E_n^{*,*}\ra E_n^{*+n,*+1-n}$. One expresses this fact by saying that $d_n$ is {\em a~differential of bidegree $(n,-n+1)$.}

Each differential $d_n$ is in a way derived from the original differential $d$, and furthermore, $E_{n+1}^{*,*}$ is the cohomology tableau of $E_n^{*,*}$ in the appropriate sense. The idea is then to compute the tableaux $E_n^{*,*}$ one by one, until they stabilize. The stabilized tableau is usually called $E_\infty^{*,*}$.

We would like to alert the reader at this point that even after the tableau $E_\infty^{*,*}$ was computed, it can still require additional work to determine the cohomology groups of the original cochain complex. Surely, if $\car$ is a~field, the situation is easy. Namely one has
        \[ H^d(\cc)=\bigoplus_{p+q=d}E_\infty^{p,q}.
\]
However, if $\car$ is an arbitrary ring (for example $\car=\dz$), then one may need to solve a~number of extension problems before obtaining the final answer. This has to do with the fact, that in a~short exact sequence of $\car$-modules
        \[0\lra A\stackrel{\alpha}{\lra} B \stackrel{\beta}{\lra}C\lra 0,
\]
$B$ does not necessarily split as a direct sum of the submodule $A$ and the quotient module $C$. This is not even true for $\car=\dz$, a~classical example is to take $A=B=\dz$, $C=\zz$, to take $\alpha:x\mapsto 2x$ to be the doubling map (injective), and to take $\beta:x\mapsto x$~$\mod 2$ to be the parity map (surjective).

Let us now describe more precisely how the tableaux $E_n^{*,*}$ and the differentials $d_n$ are constructed. As auxiliary modules, set
\[ 
Z_n^{p,q}:=C_p^{p+q}\cap d^{-1}(C_{p+n}^{p+q+1}),
\]
where $d^{-1}$ denotes the inverse of the differential $d$, i.e., $Z_n^{p,q}$ consists of all elements of $C_p^{p+q}$ whose boundary is in $C_{p+n}^{p+q+1}$; and set
\[ 
B_n^{p,q}:=C_p^{p+q}\cap d\,(C_{p-n}^{p+q-1}),
\]
i.e., $B_n^{p,q}$ consists of all elements of $C_p^{p+q}$ which constitute the image of $d$ from $C_{p-n}^{p+q-1}$. These are the settings for $n\geq 0$. Finally, for $n=-1$, we use the following convention:
\[Z_{-1}^{p,q}:=C_p^{p+q},\qquad \text{ and }\qquad 
B_{-1}^{p,q}:=d(C_{p+1}^{p+q-1}).
\]
With these notations, we set
\begin{equation} \label{eq:enpq}
 E_n^{p,q}:=Z_n^{p,q}/(Z_{n-1}^{p+1,q-1}+B_{n-1}^{p,q}),
\end{equation}
for all $0\leq n\leq\infty$.

It is an~easy check, which we leave to the reader, that $d(Z_n^{p,q})\subseteq Z_n^{p+n,q-n+1}$, and that $d(Z_{n-1}^{p+1,q-1}+B_{n-1}^{p,q})\subseteq Z_{n-1}^{p+n+1,q-n}+B_{n-1}^{p+n,q-n+1}$. It follows that, via the quotient maps, the differential $d$ induces a~map from $E_n^{p,q}$ to $E_n^{p+n,q-n+1}$, which we choose to call $d_n^{p,q}$ (or just $d_n$, if it is clear what the coefficients $p$ and $q$ are).

One can view the tableau $(E_n^{*,*},d_n)$ as a collection of (nearly) diagonal cochain complexes. This allows one to compute the cohomology groups, just like for the usual cochain complexes, by setting 
\[
H^{p,q}(E_n^{*,*},d_n)=\ker(E_n^{p,q}\stackrel{d_n}{\lra}E_n^{p+n,q-n+1})/ \im(E_n^{p-n,q+n-1}\stackrel{d_n}{\lra}E_n^{p,q}).
\]
Now we can make the sense in which $E_{n+1}^{*,*}$ is the cohomology tableau of $E_n^{*,*}$ precise:
\begin{equation} \label{eq:epqhom}
E_{n+1}^{p,q}=H^{p,q}(E_n^{*,*},d_n).
\end{equation}
Please note, that the equation~\eqref{eq:epqhom} is not trivial, and needs a~proof. It can be deduced directly from the equation~\eqref{eq:enpq}, 
see e.g., \cite{McC}. 

Let us start with unwinding these definitions for $n=0$. It follows from our conventions for $n=-1$, that 
\begin{equation} \label{eq:e0pq}
 E_0^{p,q}=(C_p^{p+q}\cap d^{-1}(C_p^{p+q+1})) /(C_{p+1}^{p+q}+d(C_{p+1}^{p+q-1}))=C_p^{p+q}/C_{p+1}^{p+q}.
\end{equation}

Furthermore, the differential $d:C_p^{p+q}\ra C_p^{p+q+1}$ induces the differential $d_0:E_0^{p,q}\ra E_0^{p,q+1}$, which is nothing else but the differential of the relative cochain complex $(\cc_p,\cc_{p+1})$. By the equation \eqref{eq:epqhom}, this yields
\begin{equation} \label{eq:e1pq}
 E_1^{p,q}=H^{p+q}(\cc_p,\cc_{p+1}).
\end{equation}
 
Moreover, one can show that $d_1^{p,q}:E_1^{p,q}\lra E_1^{p+1,q}$ is the connecting homomorphism $\partial:H^{p+q}(\cc_p,\cc_{p+1})\lra H^{p+q+1}(\cc_{p+1},\cc_{p+2})$ in the long exact sequence of the triple $(\cc_p,\cc_{p+1},\cc_{p+2})$.

Unless some additional specific information is available, it is hard to say what happens in the tableaux for $n\geq 2$. The important thing is that
with the setup above, the spectral sequence runs its course and eventually
converges (modulo the extension difficulties outlined above) to the cohomology groups of the original cochain complex.

\section{The standard spectral sequence converging to $H^*(\thomp(T,G))$.}

\subsection{Filtration induced by the support map.} \label{ssect:supp} 


\nin Let $T$ and $G$ be two graphs, and assume $T$ is finite.
As mentioned above, there is a~simplicial map $\supp:\thomp(T,G)\ra\Delta^{V(T)}$. Consider the Serre filtration of the cellular cochain complex $C^*(\thomp(T,G);\car)$ associated with this map.

 We order the vertices of $T$ and of $G$, and then observe that the vertices of $\thomp(T,G)$ are indexed with pairs $(x,y)$, where $x\in V(T)$, $y\in V(G)$, such that if $x$ is looped, then so is~$y$. Let us internally order these pairs lexicographically: $(x_1,y_1)\prec(x_2,y_2)$ if either $x_1<x_2$, or $x_1=x_2$ and $y_1<y_2$. Orient each simplex of $\thomp(T,G)$ according to this order on the vertices. We call this orientation {\it standard}, and call the oriented simplex~$\eta_+$. One can think of this simplex as a~chain in the corresponding chain complex; we denote the dual cochain with~$\eta^*_+$.

We can explicitly describe the considered filtration. Define the subcomplexes $F^p=F^pC^*(\thomp(T,G);\car)$ of $C^*(\thomp(T,G);\car)$ as follows:
\[ F^p: \dots\stackrel{\bo^{q-1}}\lra F^{p,q}\stackrel{\bo^{q}}\lra
F^{p,q+1}\stackrel{\bo^{q+1}}\lra\dots,\]
where 
\[ F^{p,q}=F^p C^q(\thomp(T,G);\car)=\car[\eta^*_+\,|\,
\eta_+\in\thomp^{(q)}(T,G),|\supp\eta|\geq p+1],
\]
$\bo^*$ is the restriction of the differential in
$C^*(\thomp(T,G);\car)$, and $\thomp^{(q)}(T,G)$ denotes the $q$-th skeleton of $\thomp(T,G)$. Phrased verbally: $F^{p,q}$ is generated by all elementary cochains corresponding to $q$-dimensional cells, which are supported in at least $p+1$ vertices of~$T$. Note, that this restriction defines a~filtration, since the differential does not decrease the cardinality of the support set.

 We have
\[
C^q(\thomp(T,G);\car)=F^{0,q}\supseteq F^{1,q}\supseteq\dots\supseteq
F^{|V(T)|-1,q}\supseteq F^{|V(T)|,q}=0,
\]
which is the Serre filtration associated to the support map.

\subsection{The 0th and the 1st tableaux.} \label{ssect:01tabl} 


\nin By writing the brackets $[-]$ after the name of a~cochain complex, we shall mean the index shifting (to the left), that is for the cochain complex $\cc=(C^*,d^*)$, the cochain complex $\cc[s]=(C^*[s],d^*)$ is defined by $C^i[s]:=C^{i+s}$; note that we choose not to change the sign of the differential.

\begin{prop} {\rm(\cite[Proposition 3.4]{BK03c}).}

\nin For any $p$,
\begin{equation}\label{eq:e0tab}
F^p/F^{p+1}=\bigoplus_{\begin{subarray}{c}{S\subseteq V(T)}\\{|S|=p+1}
\end{subarray}}
C^*(\thom(T[S],G);\car)[-p].
\end{equation}
Hence, the 0th tableau of the spectral sequence associated to
the cochain complex filtration $F^*$ is given by
\begin{equation}
  \label{eq:E0}
E_0^{p,q}=C^{p+q}(F^p,F^{p+1})=\bigoplus_{\begin{subarray}{c}
{S\subseteq V(T)}\\{|S|=p+1}\end{subarray}}
C^q(\thom(T[S],G);\car).
\end{equation}
\end{prop}

Furthermore, using the equation~\eqref{eq:e1pq}, we obtain the description of the first tableau as well.

\begin{equation} \label{eq:E1}
E_1^{p,q}=H^{p+q}(F^p,F^{p+1})=\bigoplus_{\begin{subarray}{c}
{S\subseteq V(T)}\\{|S|=p+1}\end{subarray}}H^q(\thom(T[S],G);\car).
\end{equation}

\subsection{The first differential.} \label{ssect:1diff} 

 
\nin Set now $\car=\zz$. According to the formula \eqref{eq:E1} the first differential $d_1^{p,q}:E_1^{p,q}\lra E_1^{p+1,q}$ can be viewed as a~map
\[\bigoplus_{\begin{subarray}{c}{S\subseteq V(T)}\\{|S|=p+1}\end{subarray}} H^q(\thom(T[S],G);\zz)\lra \bigoplus_{\begin{subarray}{c}
{S\subseteq V(T)}\\{|S|=p+2}\end{subarray}}H^q(\thom(T[S],G);\zz).
\]
It is possible to describe this map explicitly. 

For $S_2\subseteq S_1\subseteq V(T)$, let $i[S_1,S_2]:T[S_2] \hookrightarrow T[S_1]$ be the inclusion graph homomorphism. Since $\thom(-,G)$ is a~contravariant functor, we have an~induced map $i_G[S_1,S_2]:\thom(T[S_1],G)\ra\thom(T[S_2],G)$, and hence, an~induced map on the cohomology groups 
\[
i_G^*[S_1,S_2]:H^*(\thom(T[S_2],G);\car)\ra H^*(\thom(T[S_1],G);\car).
\]
Let $\sigma\in H^q(\thom(T[S],G);\zz)$, for some $q$, and some $S\subseteq V(G)$. The value of the first differential on $\sigma$ is given by

\begin{equation} \label{eq:d1pq}
d_1^{p,q}(\sigma)=\sum_{x\in V(T)\sm S} i_G^*[S\cup\{x\},S](\sigma).
\end{equation}

In the case of integer coefficients, $\car=\dz$, one needs more work to derive the formula for $d_1^{p,q}(\sigma)$ analogous to \eqref{eq:d1pq}, since, additionally, the signs have to be taken into consideration.

\lecture{The proof of the Lov\'asz Conjecture.}

\section{Formulation of the conjecture and sketch of the proof.}

\subsection{Formulation and motivation of the Lov\'asz Conjecture.} \label{ssect:lconj} 


\nin As mentioned in the Section~\ref{sect4}, the Lov\'asz Theorem~\ref{thm:lothm}, and the fact that the neighbourhood complex $\cn(G)$ is simply homotopy equivalent to $\bip(G)=\thom(K_2,G)$, are suggesting that $\thom$-complexes in general would provide the right context of formulating and proving further topological obstructions to graph colorings.

Up to now, the most important extension of the original Lov\'asz Theorem has been the one, where the edge $K_2$ is replaced with an~odd cycle $C_{2r+1}$.

\begin{thm} {\bf (Lov\'asz Conjecture).} \label{thm:loconj} $\,$

\nin Let $G$ be a graph, such that the complex $\thom(C_{2r+1},G)$ is $k$-connected for some $r,k\in\dz$, $r\geq 1$, $k\geq -1$, then $\chi(G)\geq k+4$.
\end{thm}

Lov\'asz Conjecture has been proved in \cite{BK03b}, for this reason we have stated it here directly as a~theorem.

\begin{rem}
It follows from Theorem~\ref{thm3.3k4}, that, once the Lov\'asz Conjecture has been proved, the statement will remain true if $C_{2r+1}$ is replaced by any graph $T$, such that $T$ can be reduced to $C_{2r+1}$ by a~sequence of folds.
\end{rem}

We formulate here a~strengthening of the original conjecture.
\begin{conj} \label{conj:myconj1}
Let $G$ be a graph, such that $\sw^k(\thom(C_{2r+1},G))\neq 0$, for some $r,k\in\dz$, $r\geq 1$, $k\geq -1$, then $\chi(G)\geq k+3$.
\end{conj}

\subsection{The winding number and the proof of the case $k=0$.} \label{ssect:wind} 


\nin  The case $k=0$ of the Lov\'asz Conjecture can be settled with relatively little machinery, that is why we choose to include for it a~separate argument. This is the "toy version" which illustrates our general methods, and we prove the more general Conjecture~\ref{conj:myconj1}.

To any continuous map $\varphi:S^1\ra S^1$ one can associate an integer $\wind(\varphi)$, called the {\em winding number} of $\varphi$. Intuitively, the absolute value of $\wind(\varphi)$ measures, as its name suggests, the number of times $\varphi$ wraps the source circle around the target circle, whereas the sign of $\wind(\varphi)$ registers whether the orientation has been changed or not. The usual way to define $\wind(\varphi)$ formally is to notice that $\varphi$ induces a~group homomorphism $\varphi^*:H^1(S^1;\dz)\ra H^1(S^1;\dz)$. Any group homomorphism from $\dz$ to itself is uniquely determined by the image of~1. This image is exactly the winding number.

As the proof of Theorem~\ref{thm:swlb} suggests, we need to analyze the complexes $\thom(C_{2r+1},K_3)$ in some detail. One can see, by direct inspection, that the connected components of $\thom(C_{2r+1},K_3)$ can be indexed by the winding numbers $\alpha$. All one needs to see is that if two homomorphisms $\varphi,\psi:C_{2r+1}\ra K_3$ have the same winding number, then there is a~sequence of edges in $\thom(C_{2r+1},K_3)$ connecting $\varphi$ with $\psi$; and this is fairly straightforward.

We notice however, that these winding numbers cannot be arbitrary.
Indeed, if the number of times $C_{2r+1}$ winds around $K_3$ is
$\alpha$, then $2r+1=3\alpha+2t$, for some nonnegative integer $t\leq
r$. Hence, $\alpha=(2r-2t+1)/3$. It follows, that $\alpha$ must be
odd, and that it cannot exceed $(2r+1)/3$. So $\alpha=\pm 1,\pm
3,\dots, \pm(2s+1)$, where $s=\lfloor(r-1)/3\rfloor$, in particular
$s\geq 0$.

Let $\varphi:\thom(C_{2r+1},K_3)\ra\{\pm 1,\pm 3,\dots,\pm(2s+1)\}$ map each point $x\in\thom(C_{2r+1},K_3)$ to the point on the real line, indexing the connected component of $x$. Clearly, $\varphi$ is a~$\zz$-map. Since $\dim(\{\pm 1,\pm 3,\dots, \pm(2s+1)\}/\zz)=0$, we have $H^1(\{\pm 1,\pm 3,\dots, \pm(2s+1)\}/\zz;\zz)=0$, and the functoriality of the Stiefel-Whitney classes implies $\sw(\thom(C_{2r+1},K_3))=0$. The Conjecture~\ref{conj:myconj1} for this case follows now from Theorem~\ref{thm:swlb}.

We have shown the Conjecture~\ref{conj:myconj1} for $k=0$ using the
Stiefel-Whitney classes, but it is equally easy to prove the Lov\'asz
Conjecture for this case directly. Indeed, following the lines of the
proof of Theorem~\ref{thm:swlb}, we see that, a~3-coloring of $G$
would induce a~$\zz$-map from $\thom(C_{2r+1},G)$ to
$\thom(C_{2r+1},K_3)$. On the other hand, the first one of these
spaces is connected, by the conjecture assumption, whereas the second
one is not, and has no connected components preserved by the
$\zz$-action. Clearly, this yields a~contradiction.

\subsection{Sketch of the proof of the Lov\'asz Conjecture.} \label{ssect:lcsketch} 


\nin Our proof of Lov\'asz Conjecture is based on two fundamental properties of the complexes $\thom(C_{2r+1},K_n)$. The first one is the following.

\begin{thm} \label{thm:swoddkn} {\rm (\cite[Theorem 2.3(b)]{BK03c}).}

\nin We have $\sw^{n-2}(\thom(C_{2r+1},K_n))=0$, for all $r\geq 1$, and odd $n$, such that 
$n\geq 3$.
\end{thm}

If $r'>r$, then there is a~$\zz$-equivariant graph homomorphism
$\varphi:C_{2r'+1}\ra C_{2r+1}$, in turn inducing a~$\zz$-map
$\varphi_{K_n}:\thom(C_{2r+1},K_n)\ra\thom(C_{2r'+1},K_n)$.  It
follows, that if Theorem \ref{thm:swoddkn} is true for $r'$,
then it is also true for $r$. Therefore, if necessary, we can assume
that $r$ is taken to be sufficiently large.

Some further details of the proof of Theorem \ref{thm:swoddkn} are
given in the Section~\ref{ssect:keven}.

To formulate the second property, consider one of the two embeddings $\iota:K_2\hra C_{2r+1}$ which maps the edge to the $\zz$-invariant edge of $C_{2r+1}$. Clearly, $\iota$ is a~$\zz$-equivariant graph homomorphism. Since $\thom(-,H)$ is a~contravariant functor, $\iota$ induces a~map of $\zz$-spaces $\iota_{K_n}:\thom(C_{2r+1},K_n)\ra \thom(K_2,K_n)$, which in turn induces a~$\dz$-algebra homomorphism
$\iota_{K_n}^*:H^*(\thom(K_2,K_n);\dz)\ra H^*(\thom(C_{2r+1},K_n);\dz)$.

\begin{thm}\label{thm:even_n} {\rm (\cite[Theorem 2.6]{BK03c}).}

\nin 
Assume $n$ is even, then $2\cdot\iota_{K_n}^*$ is a $0$-map.
\end{thm}

Some further details of the proof of Theorem \ref{thm:even_n} are
given in the Section~\ref{ssect:kodd}.

\vspace{5pt}

\noindent
{\bf Sketch of the Proof of Theorem~\ref{thm:loconj} (Lov\'asz Conjecture).}

\nin The case $k=-1$ is trivial, so take $k\geq 0$. Assume first that $k$ is even. By Corollary~\ref{crl:conn}, we have $\sw^{k+1}(\thom(C_{2r+1},G))\neq 0$. By Theorem~\ref{thm:swoddkn}, we have $\sw^{k+1}(\thom(C_{2r+1},K_{k+3}))=0$. Hence, applying Theorem~\ref{thm:swlb} for $T=C_{2r+1}$ we get $\chi(G)\geq k+4$.

Assume now that $k$ is odd, and that $\chi(G)\leq k+3$. Let $\varphi:G\ra K_{k+3}$ be a~vertex-coloring map. Combining Corollary~\ref{crl:conn}, the fact that $\thom(C_{2r+1},-)$ is a~covariant functor from loopfree graphs to $\zz$-spaces, and the map $\iota:K_2\hra C_{2r+1}$, we get the following diagram of $\zz$-spaces and $\zz$-maps:
\begin{multline*}
S_a^{k+1}\stackrel{f}{\lra}\thom(C_{2r+1},G)
\stackrel{\varphi^{C_{2r+1}}}{\lra}\thom(C_{2r+1},K_{k+3})
\stackrel{\iota_{K_{k+3}}}{\lra}\\
\stackrel{\iota_{K_{k+3}}}{\lra}\thom(K_2,K_{k+3})\cong S_a^{k+1}.
\end{multline*}
This gives a~homomorphism on the corresponding cohomology groups in
dimension $k+1$, $h^*=f^*\circ(\varphi^{C_{2r+1}})^*
\circ(\iota_{K_{k+3}})^*:\dz\ra\dz$. By Theorem~\ref{thm:even_n}
we have $2\cdot(i_{K_{k+3}})^*=0$, hence $2h^*=0$, and therefore
$h^*=0$.  It is well-known, see, e.g.,
\cite[Proposition 2B.6, p.\ 174]{Hat}, that a $\zz$-map $S_a^n\ra
S_a^n$ cannot induce a~$0$-map on the $n$th cohomology groups (in fact
it must be of odd degree). Hence, we have a~contradiction, and so
$\chi(G)\geq k+4$.  \qed

\vspace{5pt}

As the reader may have already noticed, we are actually proving a~sharper statement than the original Lov\'asz Conjecture. First of all, the condition ``$\thom(C_{2r+1},G)$ is $k$-connected'' can be replaced
by a~weaker condition ``the coindex of $\thom(C_{2r+1},G)$ is at least $k+1$''. Furthermore, for even~$k$, that condition can be weakened even further to ``$\sw^{k+1}(\thom(C_{2r+1},G))\neq 0$'', i.e., the stronger Conjecture~\ref{conj:myconj1} is proved.

\section{Completing the sketch for the case $k$ is odd.} \label{ssect:kodd} 

\subsection{The first spectral sequence and the independence complexes of cycles.}
 

\nin The main technical tool is to consider the spectral sequence 
associated to the Serre filtration induced by the support map
$\supp:\thomp(C_{2r+1},K_n)\ra\Delta^{[2r+1]}$. As we already
mentioned in the subsection~\ref{ssect:specgen}, the spectral sequence
converges to the cohomology groups of $\thomp(C_{2r+1},K_n)$. As it
happens, the complex $\thomp(C_{2r+1},K_n)$ is much easier to
understand than the complex $\thom(C_{2r+1},K_n)$.

To start with, for $n=1$, we are simply dealing with the independence
complexes of graphs: $\thomp(G,K_1)=\ind(G)$. Fortunately, that
complex has already been well-understood for cycles.

\begin{prop} \label{pr:indcr}
{\rm (\cite[Proposition 5.2]{K2})}.

\nin For any $t\geq 2$, we have
\[\ind(C_t)\simeq\begin{cases}
S^{k-1}\vee S^{k-1},&\text{ if } t=3k;\\
S^{k-1},&\text{ if } t=3k\pm 1.
\end{cases}\]
\end{prop}
Here the degenerate case $t=2$ makes sense, if we let $C_2$ be a~graph with two vertices, connected by an~edge (or a~double edge).

\begin{figure}[hbt]
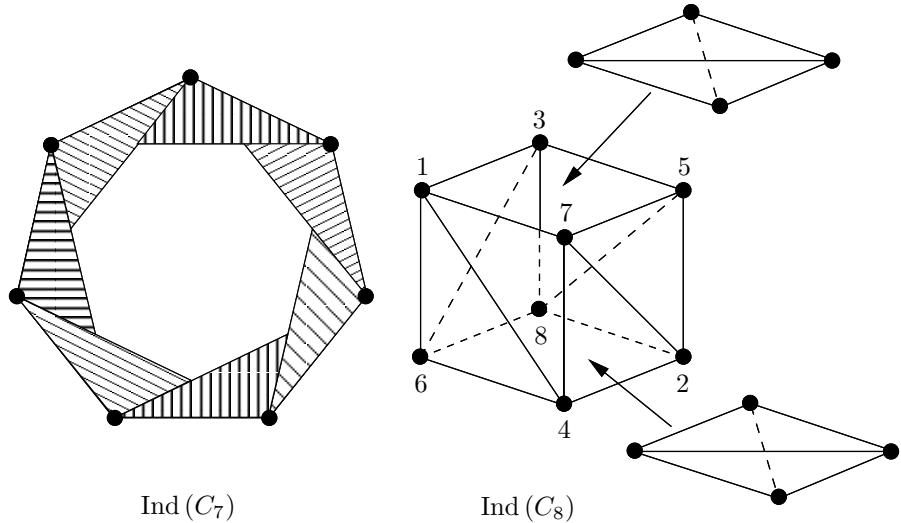

\begin{center}
  \begin{picture}(0,0)%
    \includegraphics{indc.pstex}%
  \end{picture}%
  \input{indc.pstex_t}%
 
\end{center}
\caption{Examples of independence complexes of cycles. We remark that in the right picture, the 8 triangles on the sides of the cube are filled, while the top and the bottom of the cube are filled with solid tetrahedra.}
\label{fig:indc}
\end{figure}

Now, by Corollary~\ref{crl:hompkn}, we have $\thomp(C_t,K_n)\simeq\ind(C_t)^{*n}$, hence we derive an explicit description.


\begin{crl} \label{crl:cplus}
{\rm (\cite[Corollary 4.2]{BK03c})}.

\nin For any $t\geq 2$, we have
\nin 
        \[\thomp(C_t,K_n)\simeq
        \begin{cases}
        \bigvee_{2^n\text{ copies}} S^{nk-1},&\text{ if } t=3k;\\
        S^{nk-1},&\text{ if } t=3k\pm 1.
        \end{cases}\]
\end{crl}

This is a very convenient situation for us, since we know that the spectral sequence converges to something with a~single nonzero entry.

\subsection{The analysis of the first tableau.}

Next, we look at what the first tableau of this spectral sequence is. The general formula~\eqref{eq:E1} says that the only possibly nonzero entries will be in the columns numbered $0,1,\dots,2r$. Furthermore, the entries in column number $p$ are nothing but the direct sum of the cohomology groups of induced subgraphs with $p+1$ vertices. 

For $p=2r$ this simply means that this column consists of the cohomology groups of the desired space $\thom(C_{2r+1},K_n)$ itself. For $p=0,\dots,2r-1$ we get the cohomology groups of proper induced subgraphs. Fortunately, the proper induced subgraphs of a~cycle are very simple: they are disjoint unions of isolated vertices and of strings. 

\begin{figure}[hbt]
\begin{center}
  \begin{picture}(0,0)%
    \includegraphics{1e1.pstex}%
  \end{picture}%
  \input{1e1.pstex_t}%
 
\end{center}
\caption{The $E_1^{*,*}$-tableau, for 
$E_1^{p,q}\Rightarrow H^{p+q}(\thomp(C_{2r+1},K_n);\dz)$.}
\label{fig:1e1}
\end{figure}

We recall now the formula~\eqref{eq:coprod}, which, in this particular case, says that the summands for the entries in the first tableau come from the direct products of $\thom(K_1,K_n)$ and of $\thom(L_m,K_n)$. The first one of these complexes is contractible, whereas, $L_m$ can be folded to an edge, hence, by Corollary~\ref{crl:tree}, $\thom(L_m,K_n)$ is homotopy equivalent to~$S^{n-2}$.

Since the direct products of $(n-2)$-dimensional spheres may only have nontrivial cohomology groups in dimensions which are multiples of $n-2$, we can conclude that the only possibly nontrivial entries of $E_1^{*,*}$ are in rows indexed $t(n-2)$, and in the last column. See Figure~\ref{fig:1e1} for the schematic summary of these findings; on this figure, the shaded area covers all possibly nontrivial entries of the first tableau.
 
\subsection{The analysis of the second tableau.}

Next, we need to understand what happens in every row once we pass to the second tableau. It is probably possible to perform a~complete computation. However, this is a~rather tedious task, which is unnecessary if we only care about what happens to the entries $(2r,n-2)$ and $(2r,n-3)$. Instead, we satisfy ourselves with deriving some partial information about $E_2^{*,*}$. 

The idea is to introduce some combinatorial encoding for the generators of the entries of the first tableau, and then understand the values of $d_1$. The first task is not difficult, since the generators of the cohomology groups of direct products of spheres of the same dimension, can be labeled by the subsets of the set of the spheres. 

In our case, the spheres correspond to "arcs" on the cycle, and so we can label the generators with induced subgraphs of $C_{2r+1}$, with certain set of arcs being marked. One can then filter these entries and employ another spectral sequence to compute the cohomology groups with respect to the differential~$d_1$. We refer the reader to \cite[Lemma 4.8]{BK03c} for details. 

\begin{figure}[hbt]
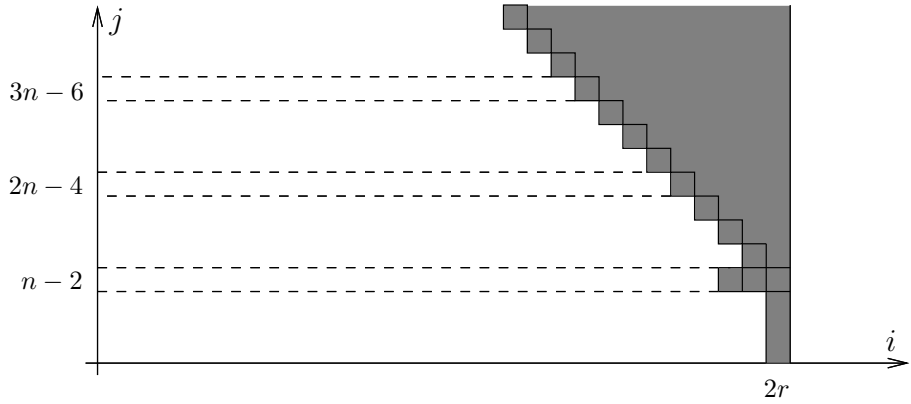

\begin{center}
  \begin{picture}(0,0)%
    \includegraphics{3e1.pstex}%
  \end{picture}%
  \input{3e1.pstex_t}%
 
\end{center}
\caption{The possibly nonzero entries in $E_2^{*,*}$-tableau,
 for $E_2^{p,q}\Rightarrow H^{p+q}(\thomp(C_{2r+1},K_n);\dz)$.}
\label{fig:3e1}
\end{figure}

The main outcome is that the possibly nontrivial entries in $E_2^{*,*}$ will be indexed by such collections of arcs, that the gaps between neighboring arcs will not exceed~2. Clearly, if we are dealing with the row $t(n-2)$, our induced subgraph of $C_{2r+1}$ cannot have fewer than $(2r+1)-2t$ vertices, since otherwise one of the gaps would be too large. Almost always this ensures that the entries of $E_2^{*,*}$ outside of the shaded area on Figure~\ref{fig:3e1} are equal to~0. There are two exceptional cases: $(n,t)=(5,2)$, and $n=4$. These cases can then be computed "by hand", using rather specific observations, see~\cite[Subsections 4.6 and 4.7]{BK03c}.
 
\subsection{The conclusion.}

After a~detailed analysis of the entries $E_2^{2r-2,n-2}$ and $E_2^{2r-1,n-2}$ we derive partial information about the cohomology groups (with integer, as well as with $\zz$ coefficients) of $\thom(C_{2r+1},K_n)$, which is summarized in the Table~\ref{tab:homol}.

\begin{table}[hbt]
\[
\begin{array}{l|c|c|c|} 
(n,r)&\,\,\,\,\,\,R\,\,\,\,\,\,& \,\,\,\,\,\,H^{n-2}\,\,\,\,\,\, &
\,\,\, \,\,\,H^{n-3}\,\,\, \,\,\,\\ \hline
2\not{|}\,\,n, n\geq 5, (n,r)\neq (5,3)\,\,\, \,\,\, & \dz &\dz&\dz\\ \hline
(n,r)=(5,3)& \dz  &\dz^2&\dz\\ \hline
2\,\,|\,\,n, n\geq 6, \text{ or } &&&\\
n=4, r\leq 3 & \dz&\zz&0\\ \hline
n=4, r\geq 4& \dz  &\dz\oplus\zz&0\\ \hline
n\geq 5, (n,r)\neq (5,3),\text{ or } &&&\\
n=4, r\leq 3& \zz &\zz&\zz\\ \hline 
(n,r)=(5,3), \text{ or } &&&\\
n=4, r\geq 4 & \zz  &{\mathbb Z}_2^2&\zz\\ \hline
\end{array}
\]
\caption{$\,$}
\label{tab:homol}
\end{table}

We remark here that the results presented in the Table~\ref{tab:homol} have been somewhat strengthened recently.
 
\begin{thm} \label{thm:ck2} \rm{(\cite[Corollary 4.6]{CK2}).}

\nin For arbitrary integers $r,n\geq 3$, the complex $\thom(C_r,K_n)$ is $(n-4)$-connected. 
\end{thm}

Let us now return to Theorem~\ref{thm:even_n}. From the Table~\ref{tab:homol}, we see that, in most of the cases, $2\cdot\iota_{K_n}^*$ is a $0$-map for a~prosaic reason: the target group $H^*(\thom(C_{2r+1},K_n);\dz)$ is isomorphic to $\zz$. The only exception is the case $n=4$, $r\geq 4$. The validity of the statement of Theorem~\ref{thm:even_n} in this special case can be verified by the direct analysis of the map $d_1:E_1^{2r-1,2}\ra E_1^{2r,2}$, see~\cite[Subsection~4.8]{BK03c} for details.

\section{Completing the sketch for the case $k$ is even.} \label{ssect:keven} 

\subsection{Topology of the quotient space $\thomp(C_{2r+1},K_n)/\zz$.}


\nin To analyze this case we need to extend some of the results of the Section~\ref{ssect:kodd}. As a~general guideline for this subsection, we would like to understand the action of $\zz$ on $\thom(C_{2r+1},K_n)$, $\thomp(C_{2r+1},K_n)$, and on their respective cohomology groups somewhat better. 

To start with, consider the $\zz$-action on $\thomp(C_{2r+1},K_n)$. Fortunately, despite of the fact, that this action is not free, it turns out to be possible to describe the quotient space rather explicitly.

\begin{prop} \label{prop:hompz2} {\rm (\cite[Proposition 4.4]{BK03c})}.

\nin    For any $r\geq 1$, we have
\[\thomp(C_{2r+1},K_n)/\zz\simeq
        \begin{cases}
        \bigvee_{2^{n-1}\text{ copies}} S^{nk-1},&\text{ if } 2r+1=3k;\\
        S^{kn/2-1}*\rp^{kn/2-1},&\text{ if } 2r+1=3k\pm 1.
        \end{cases}
\]
\end{prop}

Simple dimension inequalities yield the following corollary.

\begin{crl} \label{crl:hompz2} {\rm (\cite[Corollary 4.5]{BK03c})}.

\nin    
$\wti H^i(\thomp(C_{2r+1},K_n)/\zz)=0$ for $r\geq 2$, $n\geq 5$, and
$i\leq n+r-2$. Except for the case $r=3$.
\end{crl}

Furthermore, again by the detailed analysis of the differentials in our spectral sequence, but this time, with the $\zz$-action in mind, one can prove the following statement.

\begin{prop} \label{prop:z2act} {\rm (\cite[Corollary 4.15]{BK03c})}.

\nin    Let $n$ be odd, $n\geq 3$, $r\geq 2$, and assume $(n,r)\neq (5,3)$. Then, $\zz$ acts trivially on $H^{n-2}(\thom(C_{2r+1},K_n);\dz)$, and, it acts as a~multiplication by~$-1$, on $H^{n-3}(\thom(C_{2r+1},K_n);\dz)\simeq\dz$.
\end{prop}

We notice at this point that the support map  $\supp:\thomp(C_{2r+1},K_n)\ra\Delta^{[2r+1]}$ is $\zz$-equivariant and hence it induces the quotient map  $\supp/\zz:\thomp(C_{2r+1},K_n)/\zz\ra\Delta^{[2r+1]}/\zz$. In order to get simplicial structure on $\Delta^{[2r+1]}/\zz$, we subdivide $\Delta^{[2r+1]}$ in a~minimal way, so that every simplex preserved by $\zz$-action is fixed by this action pointwise. 

One can think of this new subdivision as the one obtained by representing simplex $\Delta^{[2r+1]}$ as a~topological join of one point and $r$ intervals: $\{c\}*[a_1,b_1]*\dots*[a_r,b_r]$, inserting an~extra vertex $c_i$ into the middle of each of the $[a_i,b_i]$, and then taking the join of $\{c\}$ and the subdivided intervals. We denote the obtained abstract simplicial complex by $\tilde\Delta^{[2r+1]}$.

The $\zz$-quotient of this simplicial structure gives one on $\Delta^{[2r+1]}/\zz$, and we can consider the Serre filtration on $\thomp(C_{2r+1},K_n)/\zz$ associated with the map $\supp/\zz$. 

\subsection{The second spectral sequence.}

Consider now the spectral sequence associated to this filtration, with the coefficients in $\zz$ instead of $\dz$. As before, this time by Proposition~\ref{prop:z2act}, we know precisely what this spectral sequence converges to. The formulae \eqref{eq:e0tab}, \eqref{eq:E0}, and \eqref{eq:E1} can be generalized as well, but before we do that we need some additional terminology.

First, we denote the set of the vertices which were added in the
subdivision by $\cc=\{c,c_1,\dots,c_r\}$. Further, for an~arbitrary
simplex $\tilde\sigma\in\tilde\Delta^{[2r+1]}$, we define its
support simplex $\vartheta(\tilde\sigma)\in\Delta^{[2r+1]}$ by replacing
every $c_i$ in $\ti\sigma$ by $\{a_i,b_i\}$, i.e.,
\[\vartheta(\tilde\sigma)=(\tilde\sigma\sm\{c_1,\dots,c_r\})\cup\bigcup_{c_i\in\ti\sigma}\{a_i,b_i\}.
\]

We can now state the analog of the formula~\eqref{eq:E1} for the spectral sequence of the quotient. The analogs of formulae \eqref{eq:e0tab}, and \eqref{eq:E0}, are straightforward, and are omitted for the sake of space, see~\cite[Section~6]{BK03c} for further details.

\begin{equation} \label{eq:qE1}
\begin{aligned}
E_1^{p,q}=& 
\bigoplus_\sigma H^{q-p}(\thom(C_{2r+1}[\vt(\sigma)],K_n)/\zz;\zz)\\
&\bigoplus_\tau H^{q-p}(\thom(C_{2r+1}[\vt(\tau)],K_n);\zz),
\end{aligned}
\end{equation}
where the first sum is taken over all $\sigma\subseteq\cc$, such that $|\sigma|=p+1$, and the second sum is taken over all $\zz$-orbits $\langle\tau\rangle$, such that $\tau\subseteq V(\ti\Delta^{[2r+1]})$, $|\tau|=p+1$, and $\tau\sm\cc\neq\emptyset$.

\begin{figure}[hbt]
\begin{center}
  \begin{picture}(0,0)%
    \includegraphics{2e1.pstex}%
  \end{picture}%
  \input{2e1.pstex_t}%
 
\end{center}
\caption{The $E_2^{*,*}$-tableau, $E_2^{p,q}\Rightarrow H^{p+q}(\thomp(C_{2r+1},K_n)/\zz;\zz)$.}
\label{fig:2e1}
\end{figure}

The next important piece of structure is understanding cohomology map in dimension $n-3$, which is induced by the quotient map $q:\thom(C_{2r+1},K_n)\ra\thom(C_{2r+1},K_n)/\zz$.

\pagebreak

\begin{prop} \label{prop:qn-3} {\rm \cite[Proposition 6.2]{BK03c}.}

\nin Let $n$ be odd, $n\geq 3$, $r\geq 2$, and assume $(n,r)\neq (5,3)$. Then,
\begin{equation}\label{eqquot}
q^{n-3}:H^{n-3}(\thom(C_{2r+1},K_n)/\zz;\zz)\ra
H^{n-3}(\thom(C_{2r+1},K_n);\zz),
\end{equation}
is a 0-map.
\end{prop}

The crucial ingredient of the proof is provided by Proposition~\ref{prop:z2act}, see~\cite{BK03c} for a~complete argument.

The proof in~\cite{BK03c} proceeds by deriving some partial information about the $E_2^{*,*}$-tableau of the spectral sequence under the consideration. The analysis is somewhat technical and we omit the details. Figure~\ref{fig:2e1} depicts the values of the entries which are of interest to us.

Let us make two important remarks. First, to derive the value $E_2^{r+1,n-3}=\zz$, one needs the result of Proposition~\ref{prop:qn-3}, which here ensures that the differential $d_1:E_2^{r,n-3}\ra E_2^{r+1,n-3}$ is a~$0$-map. Second, the value $E_2^{r-1,n-2}=0$ is derived under the assumption that $\sw^{n-2}(\thom(C_{2r+1},K_n))\neq 0$, which we are trying to disprove.

Finally, we may conclude from Figure~\ref{fig:2e1}, that $E_\infty^{r+1,n-3}=\zz$. This contradicts Corollary~\ref{crl:hompz2}, proving our original assumption $\sw^{n-2}(\thom(C_{2r+1},K_n))\neq 0$ to be wrong.

\lecture{Summary and outlook.}

\section{Homotopy tests, $\zz$-tests, and families of test graphs.}

\subsection{Homotopy test graphs.} 
\label{ss:htest} 


\nin Returning to our ideology of test graphs, it appears natural to give the following definition.

\begin{df} \label{df:htest}
A graph $T$ is called a {\bf homotopy test graph}, if, for an arbitrary graph $G$, the following equation is satisfied
\begin{equation} \label{eq:htest}
\chi(G)>\chi(T)+\conn\,\thom(T,G). 
\end{equation}
\end{df}

Using the terminology of Definition~\ref{df:htest}, Theorems~\ref{thm:cothm} and~\ref{thm:loconj} can be interpreted as saying that the complete graphs and the odd cycles are homotopy test graphs. Furthermore, it follows from Theorem~\ref{thm3.3k4}(1) that the class of homotopy test graphs is closed under the equivalence relation given by the folds and by their reverses.

\begin{figure}[hbt]
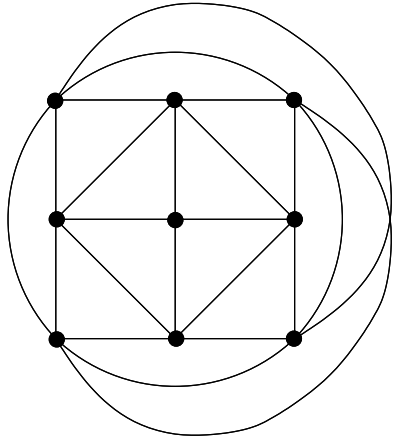

\begin{center}
  \begin{picture}(0,0)%
    \includegraphics{count.pstex}%
  \end{picture}%
  \input{count.pstex_t}%
 
\end{center}
\caption{The Hoory-Linial example of a graph, which is not a~homotopy test graph.}
\label{fig:count}
\end{figure}

More generally, it has been asked by Lov\'asz, \cite{Lo2}, whether every graph is a~homotopy test graph. That has been answered in the negative by Hoory \& Linial, \cite{HL}, whose example $HL$ is presented on Figure~\ref{fig:count}. Note that $\chi(HL)=5$, and set $G=K_5$. It was shown in \cite{HL} that $\thom(HL,K_5)$ is connected, hence the equation~\eqref{eq:htest} is false for these values of $G$ and~$T$.

The problem of characterizing the homotopy test graphs is a~formidable one, with many open questions left to explore, see for e.g., Conjecture~\ref{conj:biptest}.

\subsection{Stiefel-Whitney test graphs.} 
\label{ss:sw} 


\nin Switching from all spaces to $\zz$-spaces, and from homotopy to cohomology, we define a~different class of test graphs. First, recall the following standard notion of algebraic topology.

\begin{df} \label{df:height}
Let $X$ be a~$\zz$-space. The {\bf height} of $X$, denoted $\tth(X)$, is the maximal nonnegative integer $h$, such that $\sw^h(X)\neq 0$.
\end{df}

It is important to note, that if $X$ and $Y$ are two arbitrary $\zz$-spaces, and $\varphi:X\ra Y$ is an~arbitrary $\zz$-map, then, since the Stiefel-Whitney characteristic classes are functorial, we have $(\varphi/\zz)^*(\sw(Y))=\sw(X)$, which in particular implies the inequality $\tth(X)\leq\tth(Y)$.

We note that, for an~arbitrary $\zz$-space $X$, the existence of a~$\zz$-equivariant map $S^n_a\ra X$ implies $n=\tth(S^n_a)\leq\tth(X)$, whereas the existence of a~$\zz$-equivariant map $X\ra S^m_a$ implies $m=\tth(S^m_a)\geq\tth(X)$. This can be best summarized with the inequality
\[\coind(X)\leq\tth(X)\leq\ind(X).\]

Let us now return to graphs.

\begin{df} \label{df:swtest}
Let $T$ be a~graph with a~$\zz$-action which flips an~edge. Then, $T$ is called {\bf Stiefel-Whitney $n$-test graph}, if we have $\tth(\thom(T,K_n))=n-\chi(T)$. Furthermore, $T$ is called {\bf Stiefel-Whitney test graph} if it is Stiefel-Whitney $n$-test graph
for any integer $n\geq\chi(T)$.
\end{df}

A direct application of Theorem~\ref{thm:swlb} yields the next corollary, which also serves as an~explanation for our terminology.

\begin{crl} \label{crl:swtest}
Assume $T$ is a~Stiefel-Whitney test graph, then, for an~arbitrary graph $G$, we have 
\begin{equation} \label{eq:swtest}
\chi(G)\geq\chi(T)+\tth(\thom(T,G)). 
\end{equation}
\end{crl}

Note, that by Corollary~\ref{crl:conn}, we have $\tth(X)\geq\conn\, X+1$, for an arbitrary $\zz$-space~$X$. Therefore, comparing equations~\eqref{eq:htest} and~\eqref{eq:swtest}, we see that if a~graph $T$ is a~Stiefel-Whitney test graph, then, it is also a~homotopy test graph.

Let us stress again that, in analogy to the fact that the {\it height}
is defined for $\zz$-spaces, the term {\it Stiefel-Whitney test graph}
actually refers to a~pair $(T,\gamma)$, where $T$ is a~graph, and
$\gamma$ is an~involution of $T$, which flips an edge. The following
question arises naturally in this context.

\vskip5pt

\nin {\bf Question.} {\it Does there exist a~graph $T$ having two different 
involutions, $\gamma_1$ and $\gamma_2$, such that $(T,\gamma_1)$ is
a~Stiefel-Whitney test graph, whereas $(T,\gamma_2)$ is not?}

\vskip5pt

It would be rather surprising, if the answer to this question turned
out to be positive.

Next, we describe an important extension property of the class of
Stiefel-Whitney test graphs.

\begin{prop} \label{pr:factor}
Let $T$ be an~arbitrary graph, and let $A$ and $B$ be Stiefel-Whitney test graphs, such that $\chi(T)=\chi(A)=\chi(B)$. Assume further that there exist $\zz$-equivariant graph homomorphisms $\varphi:A\ra T$ and $\psi:T\ra B$. Then, $T$ is also a~Stiefel-Whitney test graph.
\end{prop} 

\pr Let $n$ be an arbitrary positive integer. By the functoriality of Stiefel-Whitney characteristic classes, we have 
        \[\tth(\thom(A,K_n))\leq\tth(\thom(T,K_n))\leq\tth(\thom(B,K_n)).
\]
 Hence $n-\chi(A)\leq\tth(\thom(T,K_n))\leq n-\chi(B)$, which, by the assumptions of the proposition, implies $\tth(\thom(T,K_n))=n-\chi(T)$.
\qed

\vskip5pt

The next corollary describes a~simple, but instructive example of the situation in Proposition~\ref{pr:factor}.

\begin{crl} \label{crl:z2bip}
Any connected bipartite graph $T$ with a $\zz$-action which flips an~edge is a~Stiefel-Whitney test graph.
\end{crl}

\begin{figure}[hbt]
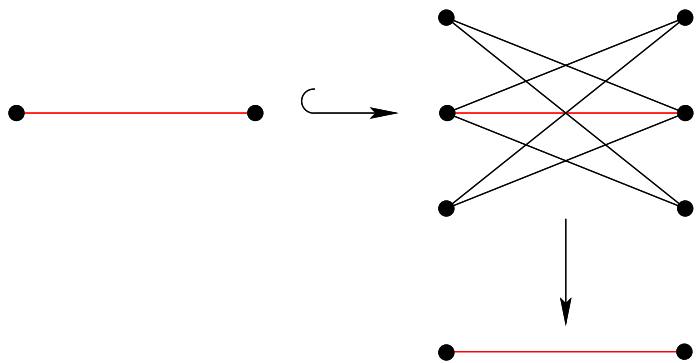

\begin{center}
  \begin{picture}(0,0)%
    \includegraphics{k2tk2.pstex}%
  \end{picture}%
  \input{k2tk2.pstex_t}%
 
\end{center}
\caption{$\zz$-invariant factoring of an edge through a bipartite graph.}
\label{fig:k2tk2}
\end{figure}

Indeed, we have $\zz$-equivariant graph homomorphisms $K_2\hra T\ra
K_2$, where the first one is the inclusion of the flipped edge, and
the second one is the arbitrary coloring map, see
Figure~\ref{fig:k2tk2}. Since, by Proposition~\ref{pr:k2kn}, $K_2$ is
a~Stiefel-Whitney test graph, we conclude that $T$ is also
a~Stiefel-Whitney test graph.

In particular, any even cycle with the $\zz$-action which flips
an~edge is a Stiefel-Whitney test graph.

\vspace{5pt}
\nin {\it Summary.} The class of Stiefel-Whitney test graphs contains 
complete graphs, connected bipartite graphs (in both cases one can
take any involution which flips an edge). Furthermore, it is closed
under factorizations, as described in Proposition~\ref{pr:factor}.

By Theorem~\ref{thm:swoddkn}, the odd cycles are Stiefel-Whitney
$n$-test graphs, for odd $n\geq 3$. Conjecturally, see
Conjecture~\ref{conj:swmain}, odd cycles are Stiefel-Whitney $n$-test
graphs, for all $n\geq 3$.

\section{Conclusion and open problems.}

\nin It follows from Corollary~\ref{crl:z2bip} that any connected bipartite graph with a $\zz$-action which flips an edge is a~homotopy test graph. It seems natural to generalize this statement.

\begin{conj} \label{conj:biptest}
Every connected bipartite graph is a homotopy test graph.
\end{conj}

By what is said above, we know that a~connected bipartite graph is a~homotopy test graph if there is a~sequence of folds and their reverses, reducing it to some connected bipartite graph with a $\zz$-action which flips an~edge.

Before formulating the next conjecture, we recall that by saying that a~topological space is $(-1)$-connected, we mean that it is nonempty. Clearly, if the maximal valency of $G$ is at most $n-1$, then $G$ can be colored with $n$ colors by means of the greedy procedure. Furthermore, Babson \& Kozlov proved in \cite[Proposition 2.4]{BK03b} that if the maximal valency of $G$ is at most $n-2$, then $\thom(G,K_n)$ is $0$-connected. Generalizing this statement to higher dimension, we obtain the following conjecture.

\begin{conj} \label{conj:connh} {\rm (Babson \& Kozlov, \cite[Conjecture~2.5]{BK03b}).}\footnote{At the time of the writing of this survey, this conjecture has been proved and is now a~theorem, see~\cite{CK2}.}

\nin Let $G$ be any graph.  If the maximal valency of $G$ is equal to~$d$, then $\thom(G,K_n)$ is $k$-connected, for all integers $k\geq -1$, $n\geq d+k+2$. 
\end{conj}

Next, let us recall an important class of manifolds.

\begin{df} \label{df:stiman}
For an arbitrary positive integer $n$, the {\bf Stiefel manifold} $V_k(\rn)$ is the set of the orthonormal $k$-frames in an~$n$-dimensional Euclidean space, topologized as subspace of $(\rn)^k$.
\end{df}

Stiefel manifolds are homogeneous spaces and play an important role in the study of characteristic classes, see~\cite{MS}.

\begin{conj} \label{conj:csorba} {\rm (Csorba, \cite{Cs2}). }

\nin The complex $\thom(C_5,K_n)$ is homeomorphic to $V_2(\br^{n-1})$, for all $n\geq 1$. 
\end{conj}

The cases $n=1,2$ are tautological, as both spaces are empty. The example on the Figure~\ref{fig:c45k3} verifies the case $n=3$: $\thom(C_5,K_3)\cong S^1\coprod S^1$. Several cases, including $n=4$ have been recently verified by Csorba \& Lutz, see~\cite{CL1}.

\vskip7pt

Returning to the Stiefel-Whitney characteristic classes, we have the following hypothesis.

\begin{conj} \label{conj:swmain} {\rm (Babson \& Kozlov, \cite[Conjecture~2.5]{BK03c}).}

\nin The equation
\begin{equation} \label{eq:ccvan} 
\sw^{n-2}(\thom(C_{2r+1},K_n))=0, \text{ for all }n\geq 2
\end{equation}
is true for an arbitrary positive integer~$r$.
\end{conj}

Clearly, the case $n=2$ is obvious, since $\thom(C_{2r+1},K_2))=\emptyset$.
The Conjecture~\ref{conj:swmain} has been proved in \cite{BK03c} for $r=1$ and arbitrary $n\geq 2$, as well as for odd~$n$ and arbitrary~$r$, see here Theorem~\ref{thm:swoddkn}. For $r=2$, $n=4$, the equation~\eqref{eq:ccvan} follows from the fact that $\thom(C_5,K_4)\cong\rp^3$, and the analysis of the corresponding $\zz$-action on $\rp^3$. 

We remark here that the Conjecture~\ref{conj:swmain}, coupled with Theorem~\ref{thm:swlb}, implies the Conjecture~\ref{conj:myconj1}. Note also that, as previously remarked, for a~fixed value of $n$, if the equation \eqref{eq:ccvan} is true for $C_{2r+1}$, then it is true for any $C_{2\ti r+1}$, if $r\geq\ti r$.

\vskip7pt

We finish with another conjecture by Lov\'asz. In \cite{BrW}, Brightwell \& Winkler have shown the following result.

\begin{thm} {\rm (Brightwell \& Winkler, \cite{BrW}).}

\nin Let $G$ be an arbitrary graph. If for any graph $T$, with maximal
valency at most $d$, the graph $\thom_1(T,G)$ is connected or empty, then
$\chi(G)\geq \frac{d}{2}+2$.
\end{thm} 

Lov\'asz has suggested that this statement can be strengthened, and
that furthermore, a higher dimensional analog is true.

\begin{conj} {\rm (Lov\'asz).}

\nin Let $G$ be an arbitrary graph. If for any graph $T$, with maximal
valency at most $d$, the complex $\thom(T,G)$ is $k$-connected or empty, then $\chi(G)\geq d+k+2$.
\end{conj}

\vfill



\nin {\bf Acknowledgements.} We thank Peter Csorba, Sonja \v{C}uki\'c, 
Alexander Engstr\"om, as well as the anonymous referees, for the
helpful comments concerning the presentation in this paper.

\end{document}